\documentclass[11pt]{article}

\usepackage{amsmath}
\usepackage{graphicx}
\usepackage{grffile}
\usepackage{latexsym}
\usepackage{color}
\usepackage{amssymb}
\usepackage{amsbsy}
\usepackage{hyperref}
\usepackage[linesnumbered,ruled,vlined,commentsnumbered]{algorithm2e}
\usepackage{authblk}
\usepackage{geometry}
\usepackage{epstopdf}
\usepackage{cite}
\usepackage{booktabs}
\usepackage{array,multirow}

\newtheorem{remark}{Remark}

\newtheorem{theorem}{Theorem}

\newtheorem{lemma}{Lemma}

\def\RR{\mathbb R}
\def\e{\varepsilon}

\def\U{{\bf u}}
\def\V{{\bf v}}
\def\W{{\bf w}}
\def\Z{{\bf z}}
\def\F{{\bf f}}
\def\FF{{\bf F}}
\def\G{{\bf g}}
\def\E{{\bf E}}
\def\D{{\bf D}}

\def\be{\begin{equation}}
\def\ee{\end{equation}}
\def\bea{\begin{eqnarray}}
\def\eea{\end{eqnarray}}
\begin{document}
\title{Gradient-based Monte Carlo methods for relaxation approximations of hyperbolic conservation laws}

\author[$*$]{Giulia Bertaglia}
\author[$\dagger$]{Lorenzo Pareschi}
\author[$\ddagger$]{Russel E. Caflisch}
\affil[$*$]{\it Department of Environmental and Prevention Sciences, University of Ferrara, Italy}
\affil[$\dagger$]{\it Department of Mathematics and Computer Science, University of Ferrara, Italy}
\affil[$\ddagger$]{\it Courant Institute of Mathematical Sciences, New York University, USA}

\maketitle

\begin{abstract}
Particle methods based on evolving the spatial derivatives of the solution were originally introduced to simulate reaction-diffusion processes, inspired by vortex methods for the Navier--Stokes equations. Such methods,  referred to as gradient random walk methods, were extensively studied in the '90s and have several interesting features, such as being grid free, automatically adapting to the solution by concentrating elements where the gradient is large and significantly reducing the variance of the standard random walk approach.
In this work, we revive these ideas by showing how to generalize the approach to a larger class of partial differential equations, including hyperbolic systems of conservation laws. To achieve this goal, we first extend the classical Monte Carlo method to  relaxation approximation of systems of conservation laws, and subsequently consider a novel particle dynamics based on the spatial derivatives of the solution. The methodology, combined with asymptotic-preserving splitting discretization, yields a way to construct a new class of gradient-based Monte Carlo methods for hyperbolic systems of conservation laws.
Several results in one spatial dimension for scalar equations and systems of conservation laws show that the new methods are very promising and yield remarkable improvements compared to standard Monte Carlo approaches, either in terms of variance reduction as well as in describing the shock structure.
\end{abstract}


\tableofcontents

\section{Introduction}
Monte Carlo methods, which were first devised during the Manhattan Project in the 1940s to simulate the behavior of neutron diffusion in fissile material~\cite{Nick49, Ulam47}, have a long and storied history~\cite{Nick87}.
In the decades following their introduction, the methods were refined and extended to various fields, including physics, engineering, finance, and many others. The subsequent development of computers provided a significant boost to the spread of Monte Carlo methods in scientific computing, as they allowed for the simulation of more complex systems and the generation of larger numbers of random samples. 
Although today Monte Carlo methods are widely used in many scientific and industrial applications~\cite{MC}, their systematic design as a numerical analysis tool to solve partial differential equations (PDEs) is still rather limited to specific contexts compared to deterministic approaches such as finite-differences or finite-volumes. Prominent examples are applications involving  stochastic differential equations, like in financial applications~\cite{finance1, finance2, finance3, Caflisch98}, or when other numerical methods are not feasible due to the complexity of the problem, like in plasma physics and rarefied gas dynamics~\cite{bird,Nanbu80,CPmc,PRMC,Trazzi,dimarco1, RW05}.

On the other hand, growing attention in high-dimensional problems and in quantifying uncertainty in many fields of applied mathematics, including emerging fields such as life and social sciences, has greatly increased interest in developing efficient Monte Carlo approaches in such contexts~\cite{giles, dimauq, Hu21, Ber, PareschiToscani}. Moreover, since the Monte Carlo solution is computed by averaging independent calculations, the resulting algorithms are well suited to exploit modern parallel computing techniques.

Among the various Monte Carlo methods developed for PDEs, one approach, inspired by the so-called vortex method for the Navier-Stokes equations~\cite{Chorin73}, had considerable success in the 1990s for reaction-diffusion problems. The idea of the method is to use the spatial derivative of the solution (i.e., the gradient in multiple dimensions) as the unknown variable. This allows the statistical solution in the original variables to be reconstructed not from the histogram, as in standard random walk approaches, but directly as the distribution function of the samples of the derivative. In addition to the statistical reduction of fluctuations, such a technique, usually referred to as the Gradient Random Walk (GRW)  method, brings with it the advantages of adaptivity, since samples are taken according to the space derivatives, and a grid-free structure~\cite{Rob, SM94, SP86, GS, CR90}.

In this paper we explore the possibility of extending such ideas to a broader class of PDEs that includes, in particular, systems of conservation laws. 
The construction of Monte Carlo methods for such problems is strongly inspired by the fluid-dynamical limit of the Boltzmann equation and the corresponding Monte Carlo methods for the Euler equations~\cite{CPmc, PRMC, Trazzi, PRMC2, pullin}. Let us remark that the construction of stochastic particle methods for nonlinear hyperbolic problems has been the subject of limited research in the past and that there is no general approximation methodology (see for example~\cite{BoTa, Rob, Masc19}). A first attempt in this direction, limited to nonnegative solutions, has been proposed in~\cite{Seaid} based on appropriate relaxation approximations of such systems~\cite{JinXin, ADN, Jin2}. Indeed, such a relaxation approximation allows the Monte Carlo strategy traditionally used for kinetic equations to be adapted to generic conservation laws. We mention here that related approaches based on deterministic particles have been also proposed~\cite{Chertock, Russo1990, difrancesco2017, DFR2017, Mascagni1995}.

Following the idea of approximating the system of conservation laws by means of a semi-linear hyperbolic system with source terms we show how it is possible to formulate an approach based on statistical samples of the spatial derivatives of the solution that allows the characteristic advantages of GRW methods for reaction-diffusion equations to be generalized to such problems. In order to compare the new methods, which we will refer to as Gradient-based Monte Carlo (GBMC) methods, with a standard Monte Carlo approach, we will recall the concepts introduced in~\cite{Seaid} by extending them to the case of solutions of arbitrary sign and discussing some  variance reduction strategies. 

Here, we will limit ourselves to the one-dimensional case, which will enable us to focus on the fundamental concepts of the GBMC method and address both the scalar case and the case of systems of conservation laws, leaving to further research multidimensional extensions. We will show how the resulting methods will be grid-free in the scalar case, while in the case of systems, this may not always be possible. In all test cases considered, however, the GBMC method offers significant advantages over the traditional Monte Carlo approach, including better spatial resolution of discontinuities and a reduction in variance by several orders of magnitude.

The rest of the article is organized as follows. In Section~\ref{sec:GRW} we will recall the basic concepts of GRW methods for reaction-diffusion problems. Then in Section~\ref{sec:GBMCsc} we will discuss the design of the novel GBMC approach in the case of simple relaxation systems for scalar conservation laws. In this section we will also introduce the details of a direct Monte Carlo approach by generalizing some of the ideas presented in~\cite{Seaid}. Several numerical examples for various scalar problems illustrate the strong advantages of GBMC over standard Monte Carlo. Next, Section~\ref{sec:GBMCsy} is devoted to the extension of the methodology to the case of hyperbolic relaxation approximations to systems of conservation laws. In such cases, even if a grid free approach is related to the possibility to diagonalize the system, we show through several examples that the resulting GBMC schemes maintains most of the advantages observed in the scalar case. Some concluding remarks and further research directions are discussed in the last section. A separate appendix reports some estimates in the $L^p$ distance of the solution reconstruction error from the particles showing the improvement in accuracy of the GBMC approximation over the corresponding MC.

\section{Gradient random walk for reaction-diffusion equations}
\label{sec:GRW}
In this section we recall the basic ideas behind the design of GRW methods for reaction-diffusion problems~\cite{Mascagni1995,SM94,Rob,SP86,CR90,GS}. 
We illustrate the idea  by treating the one-dimensional scalar reaction--diffusion equation
\be
\begin{aligned}
&\frac{\partial u}{\partial t}=D\displaystyle\frac{\partial^2 u}{\partial x^2}+G(u),\quad x\in \RR,\quad D, t>0,\\
&u(x,0)=u_0(x).
\label{eq.react-diff}
\end{aligned} 
\ee
Let us assume that $u_0(x)$ is monotone and bounded from below and above. Without loss of generality, we impose $u(-\infty,t)=0$ and $u(\infty,t)=1$, for all $t\ge 0$.
We introduce the auxiliary variable $w=\partial u/\partial x$ and observe that it satisfies the equation
\be
\begin{aligned}
&\frac{\partial w}{\partial t}=D\displaystyle\frac{\partial^2 w}{\partial x^2}+G'(u)w,\quad x\in \RR,\quad D, t>0,\\
&w(x,0)=\frac{\partial u_0(x)}{\partial x}.
\label{eq.react-diff2}
\end{aligned} 
\ee
Note that, by virtue of the assumptions on $u_0(x)$, $w_0(x)$ is a probability density. 


Using $N$ particle samples located respectively in $X_1,\ldots,X_N$ (we use capital $X$ to indicate that the positions are random variables) we will have
\be
w_N(x,t)=\frac1{N}\sum_{k=1}^N \delta(x-X_k),
\ee
discretizing $w$ as a sum of $\delta$-functions. Then, since $u(x,t)=\int_{-\infty}^x w(y,t)\,dy$ is the cumulative distribution function of the random variable $X$, we get
\be
u_N(x,t)=\frac1{N}\sum_{k=1}^N H(x-X_k),
\label{eq.recL}
\ee
where $H(\cdot)$ is the (right-continuous) Heaviside step function
 \[
H(x) = 0 \,\,\, {\rm if}\,\,\, x <0,  \qquad H(x) = 1 \,\,\, {\rm if}\,\,\, x \ge 0.
 \]
Note that the boundary condition at the left is exactly satisfied, while the boundary condition at the right is satisfied on average, with fluctuations. Of course we can use a right reconstruction approach that would lead to the opposite situation (see also Remark \ref{rmk.REC}).

\begin{figure}[t]
\centering
\includegraphics[width=0.32\textwidth]{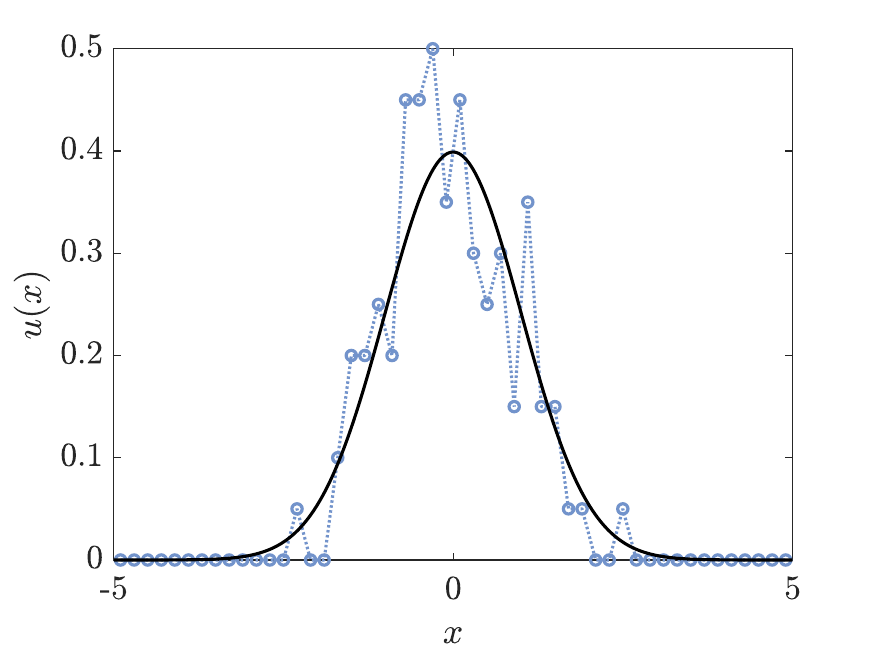}
\includegraphics[width=0.32\textwidth]{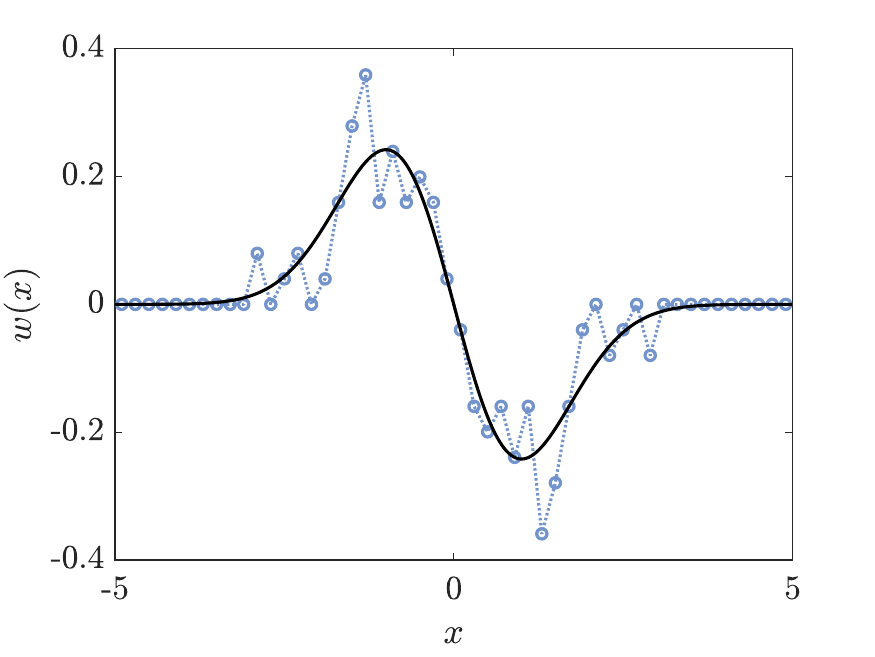}
\includegraphics[width=0.32\textwidth]{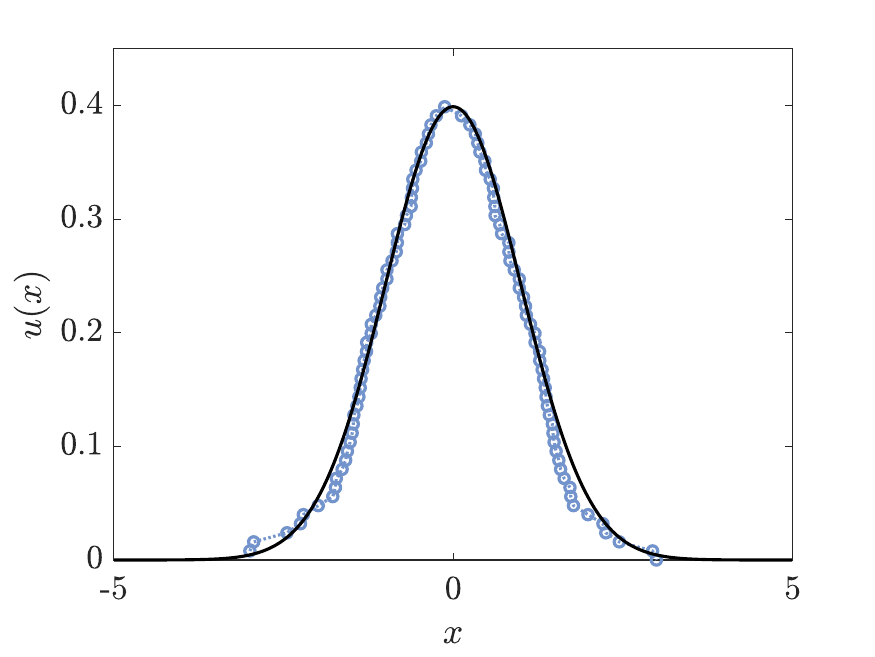}
\caption{Standard normal distribution using $N=100$ samples. Left: histogram of Monte Carlo samples with $M=50$ cells. Middle: histogram of its derivative using two symmetric families of samples with positive and negative masses and $M=50$ cells. Right: plot of the cumulative distribution of the derivative samples using left reconstruction as in \eqref{eq.recL}.
}
\label{RW_reconstruction}
\end{figure}

Notice that the same idea can be applied in case of weighted particles, namely particles having different {mass} $m_i$, thus approximating the function $w(\cdot, t)$ as 
\be\label{eqch3:vv}
w_N(x,t)=\frac1{N}\sum_{k=1}^N m_k \delta(x-X_k).
\ee
Here the mass values $m_i\in \{-m,m\}$, $m>0$ might also be negative, so we can approximate more than monotone non-negative solutions. From definition \eqref{eqch3:vv}, in this case by left reconstruction we obtain
\be
u_N(x,t)=\frac1{N}\sum_{k=1}^N m_k H(x-X_k).
\label{eqch3:cs3}
\ee

The computed solution for $u(x,t)$ contains much less fluctuation than $w(x,t)$, since all particles contribute to the solution at any point $x$, as depicted in Figure~\ref{RW_reconstruction}. In particular, due to its nature of solving the problem looking at the gradient of the state variable, the GRW method perfectly applies to problems in which sharp fronts appear and where one needs to spatially resolve jump-like solutions, because the density of the particles is large precisely where $u$ has large gradients.
Once the initial data have been discretized, the GRW method evolves the positions and masses of the particles so that $u$ satisfies \eqref{eq.react-diff}. This is done by a fractional-step iteration in a small time interval $\Delta t$ in which the diffusion term is modeled by a random walk and the reaction term is modeled through a particle-killing or replication approach with probability $\Delta t |G'(u)|$  (see~\cite{Mascagni1995,SM94,SP86,GS} for more details).

\section{Gradient-based Monte Carlo for hyperbolic problems}
\label{sec:GBMCsc}
The extension of the previous idea to hyperbolic problems and in general to nonlinear PDEs is nontrivial and, except for some specific cases~\cite{Rob}, has been poorly studied in the literature. Among other reasons, it should be pointed out that already the construction of a standard Monte Carlo method is not obvious in the case of nonlinear PDEs~\cite{Seaid, BoTa, pullin}. In the following we will consider the case of hyperbolic problems, with special reference to the relaxation approximation of scalar conservation laws. In the following we will consider the case of hyperbolic problems, with special reference to the approximation by relaxation of scalar conservation laws. Subsequently, we will discuss how to extend the method to the case of systems of conservation laws.  

The starting point is the following hyperbolic system with relaxation introduced in~\cite{JinXin} as a semi-linear approximation to scalar conservation laws 
\be
\begin{aligned}
\frac{\partial u}{\partial t} + \frac{\partial v}{\partial x} &=  0,\\
\frac{\partial v}{\partial t} + a^2 \frac{\partial u}{\partial x} &= -\frac1{\e} \left(v-F(u)\right),
\end{aligned}
\label{eqch3:jx}
\ee
where $x\in\Omega\subseteq \mathbb{R}$, supplemented with the initial conditions $u(x,0)=u_0(x)$, $v(x,0)=v_0(x)$ and suitable boundary conditions.

In the limit $\e\to 0$ from the second equation in \eqref{eqch3:jx} one formally obtains the {local equilibrium}
$v=F(u)$ and thus solutions to (\ref{eqch3:jx}) are well approximated by the scalar conservation law
\be
\frac{\partial u}{\partial t} + \frac{\partial F(u)}{\partial x} = 0.
\label{eqch3:cls}
\ee
More precisely, if we evaluate the $\e$ perturbation of the local equilibrium
\begin{equation*}
v = F(u) + \e v_1,
\end{equation*}
and substite in the second equation in \eqref{eqch3:jx}, with some algebraic manipulations we obtain
\begin{equation*}
v_1 
 = F'(u)^2\frac{\partial u}{\partial x} - a^2\frac{\partial u}{\partial x} + O(\e),
\end{equation*}
and therefore
\begin{equation*}
v = F(u) + \e\left( F'(u)^2 - a^2\right)\frac{\partial u}{\partial x},
\end{equation*}
which leads to the following convection-diffusion equation
\be
\frac{\partial u}{\partial t} + \frac{\partial F(u)}{\partial x} = \e \left[\left( a^2- F'(u)^2\right)\frac{\partial u}{\partial x}\right]_x.
\ee
The above equation is a good approximation of system \eqref{eqch3:jx} for $\e \ll1$ only if the characteristic velocity of the relaxed system is dominated by the characteristic velocity of the original system, hence if the {subcharacteristic condition} $a^2 > F'(u)^2$ is satisfied (see~\cite{CLL, JinXin}).

It is worth to underline that, for small values of $\e$, the study of  a numerical solution of system (\ref{eqch3:jx}) requires particular care even though the transport operator is linear because a standard explicit discretization will lead to time constraints of the order of $\varepsilon$. This idea was at the basis of the so-called relaxation schemes for hyperbolic systems of conservation laws~\cite{JinXin}.

\subsection{A direct Monte Carlo approach}
\label{ch4:jinxin}
First, we recall the probabilistic approach proposed in~\cite{Seaid} to construct a random particle solver for nonnegative densities that works uniformly with respect to the relaxation rate $\e$. The method is inspired by the classical Direct Simulation Monte Carlo (DSMC) method for rarefied gas dynamics~\cite{bird, Nanbu80, PRMC}. We will later discuss how to extend this approach to solutions of arbitrary sign.

Given $a>0$, let us rewrite system \eqref{eqch3:jx} introducing the diagonal (kinetic) variables
\[
f^+=\frac{{a}u+v}{2{a}}, \qquad f^-=\frac{{a}u-v}{2{a}}, \,\,
\] and the corresponding equilibrium states
\[
E^+(u)=\frac{{a}u+F(u)}{2{a}}, \qquad E^-(u)=\frac{{a}u-F(u)}{2{a}}.
\]
We obtain the diagonal form of the relaxation system 
\be
\begin{aligned}
\frac{\partial f^+}{\partial t} + {a}\frac{\partial f^+}{\partial x} &= -\displaystyle\frac1{\e} \left(f^+-E^+(u)\right),\\
\frac{\partial f^-}{\partial t} - {a} \frac{\partial f^-}{\partial x} &= -\displaystyle\frac1{\e} \left(f^--E^-(u)\right).
\end{aligned}\label{ch4:gt}
\ee
System \eqref{ch4:gt} presents several analogies with discrete-velocity models of the Boltzmann equation, in the sense that it describes a system of particles having only two speeds $\pm a$ which relax toward the local equilibrium states $E^{\pm}(u)$. The kinetic interpretation requires $E^{\pm}(u)\geq
0$ and this is guaranteed if $a \ge |F(u)|/u$ and $u > 0$. Under these assumptions,
one can apply a direct simulation Monte Carlo approach~\cite{Seaid}. 

More precisely, the solution in a small time interval $[0, \Delta t]$ is approximated by means of a fractional step procedure that solves separately the two problems characterized by the linear transport and by the relaxation term. One therefore solves, one after the other, a transport step
\be
\begin{aligned}
\frac{\partial f^+}{\partial t} + {a}\frac{\partial f^+}{\partial x} &=0,\\
\frac{\partial f^-}{\partial t} - {a} \frac{\partial f^-}{\partial x} &=0,
\end{aligned}\label{ch4:gttra}
\ee
and a relaxation process
\be
\begin{aligned}
\frac{\partial f^{+}}{\partial t} &= -\frac1{\e} \left(f^{+}-E^+(u)\right),\\
\frac{\partial f^{-}}{\partial t} &= -\frac1{\e} \left(f^{-}-E^-(u)\right).
\end{aligned}
\label{eqch3:os}
\ee
Since both steps, taken separately, can be solved exactly no further approximation besides the $O(\Delta t)$ error of the above splitting, is necessary. In fact, in a small time interval $\Delta t$ the exact solutions $\tilde
f^{\pm}$ of the free transport (\ref{ch4:gttra}) reads
\be
\begin{aligned}
\tilde{f}^+(x,\Delta t)&=f^+_0(x-a\Delta t),\\
\tilde{f}^-(x,\Delta t)&=f^-_0(x+a\Delta t),
\end{aligned}
\ee
and, setting $\tilde{u}(x,\Delta t)=\tilde{f}^+(x,\Delta t)+\tilde{f}^-(x,\Delta t)$, the exact solution of the relaxation step (\ref{eqch3:os}) yields the approximated values at time $\Delta t$
\be
\begin{aligned}
f^{+}(x,\Delta t)&=e^{-\Delta t/\e}\tilde{f}^+(x,\Delta t)+\left(1-e^{-\Delta t/\e}\right)E^+(\tilde{u}(x,\Delta t)),\\
f^{-}(x,\Delta t)&=e^{-\Delta t/\e}\tilde{f}^-(x,\Delta t)+\left(1-e^{-\Delta t/\e}\right)E^-(\tilde{u}(x,\Delta t)).
\end{aligned}
\label{ch4:rel}
\ee
The probabilistic interpretation is readily obtained by introducing the new variables 
\be
p^+(x,t)=\frac{f^{+}(x,t)}{u(x,t)},\qquad
p^-(x,t)=\frac{f^{-}(x,t)}{u(x,t)}, \label{ch4:nd}\ee with
$p^+(x,t)+p^-(x,t)=1$, to have the convex combination
\be
\begin{aligned}
p^+(x,\Delta t)&=e^{-\Delta t/\e}{p_0^+}(x-a\Delta t)+\left(1-e^{-\Delta t/\e}\right)p^+_E(x,\Delta t),\\
p^-(x,\Delta t)&=e^{-\Delta t/\e}{p_0^-}(x+a\Delta t)+\left(1-e^{-\Delta t/\e}\right)p^-_E(x,\Delta t),
\end{aligned}
\label{ch4:rel3}
\ee
where we denoted the equilibrium probabilities as
\be
\label{eq:proj}
p^+_E(x,t)=\frac{E^+(\tilde{u}(x,t))}{\tilde{u}(x,t)},\qquad
p^-_E(x,t)=\frac{E^-(\tilde{u}(x,t))}{\tilde{u}(x,t)}. 
\ee 
We note that the probabilistic interpretation holds true independently of the choice of $\Delta t/\e$. In particular, as $\varepsilon\to 0$ it reduces to the projections
\be
\label{pE}
p^+(x,\Delta t)=p^+_E(x,\Delta t),\qquad
p^-(x,\Delta t)=p^-_E(x,\Delta t),
\ee
which characterize the probabilistic approximation of the limit equation (\ref{eqch3:cls}).

A Monte Carlo method can be derived by sampling directly from the exact solutions of the operator splitting steps (\ref{ch4:rel3}). Note, however, that the probability of a velocity change in \eqref{eq:proj} depends on the mass density after the transport step $\tilde u(x,t)$.
Therefore, in order to set up a Monte Carlo method and to estimate the probability of a velocity change, we must reconstruct the mass density $u$ in a neighborhood of the particle position. 
 Given a set
of samples $X_1,\ldots,X_N$ the simplest method, which produces a piecewise constant
reconstruction, is based on evaluating the histogram of the
samples at the cell centers of a suitable grid $x_j=j\Delta x$,
\be
u_j(t)=\frac{1}{N}\sum_{k=1}^N \Phi_{\Delta x}(x_{j+1/2}-X_k),\qquad
j=\ldots,-2,-1,0,1,2,\ldots
\label{eq:rhor}
\ee
where {$\Phi_{\Delta x}(x)=1/\Delta x$} if {$|x|\leq \Delta x$}
and {$\Phi_{\Delta x}(x)=0$} elsewhere. 
Smoother versions can be obtained by changing the function approximating the Dirac delta. This corresponds to a convolution of the samples with a suitable mollifier~\cite{LP, PRMC}. 
%


We can set up a Monte Carlo method in the case of positive solutions as follows.
Given a set of samples
$(X^0_1,V^0_1), \ldots, (X^0_{N},V^0_N)$, where the particles velocities $V^0_i~\in~\{-a,a\}$ distinguish the samples of $f^+_0(x)$ from
those of $f^-_0(x)$, a new set of samples $(X_1,V_1),\ldots,(X_{N},V_N)$ is generated by Algorithm \ref{alg.MCscalar}. 
\begin{algorithm}[h!]
\label{alg.MCscalar}
\caption{Monte Carlo for $2\times 2$ hyperbolic relaxation systems}
\begin{enumerate}
\item Compute
$X_i=X_i^0+V_i^0 \Delta t,\quad i=1,\ldots,N$.
\item At the cell centers of a space grid $x_j=j\Delta x$, reconstruct $u_j(t)$
from (\ref{eq:rhor}).
\item In each space cell $j$, given a sample $(X_i,V_i^0)$
\begin{enumerate}
\item with probability $1-e^{-\Delta t/\e}$ do the following:
\begin{itemize}
\item[-] with probability $\displaystyle\frac{E^+(u_j(t))}{u_j(t)}$ set
$V_i=a$,
\item[-] otherwise set $V_i=-a$.
\end{itemize}
\item otherwise set $V_i=V_i^0$.
\end{enumerate}
\end{enumerate}
\end{algorithm}


Note that there is no need for the grid to be uniform, and the choice of the grid points can easily change at any time step. Although in principle using higher order fractional step methods would lead to an increase in the order of accuracy in $\Delta t$, in the limit $\e\to 0$ such methods, as it is well known, degenerate to first order \cite{Jin2, PRMC, PRMC2}. It is to date an open problem to construct splitting techniques of order higher than one that maintain accuracy in the limit $\e \to 0$.

Algorithm \ref{alg.MCscalar} can be improved by adopting a simple low variance technique which estimates the number of particles that should have the different velocities before their assignment. This can be done by replacing step 3 by first computing in each cell $j$ the total number of interacting particles 
\be
N^c_j(t)= {\rm SRound}\left((1-e^{-\Delta t/\e})N_j(t)\right),
\label{eq:Nc}
\ee
and then assign the velocities $a$ and $-a$ respectively to a number 
\be
N_j^+(t) = {\rm SRound}\left(\frac{E^+(u_j(t))}{u_j(t)}{N}^c_j(t)\right),\quad N_j^-(t)=N^c_j(t)-N_j^+(t),
\label{eq:varr}
\ee
of randomly chosen particles in cell $j$. An implementation is reported in Algorithm \ref{alg.VarReduction}.

%
\begin{algorithm}[h!]
\label{alg.VarReduction}
\caption{Monte Carlo variance reduction in relaxation step}
\begin{enumerate}\setcounter{enumi}{2}
\item In each space cell $j$, 
\begin{enumerate}
\item compute $N_j^+(t)$ and $N_j^-(t)$ accordingly to \eqref{eq:Nc}-\eqref{eq:varr};
\item perform a random permutation $\{k_1,\ldots,k_{N_j(t)}\}$ of the particle indices;
\begin{itemize}
\item[-] set $V_i=a$ for $i=k_1,\ldots,k_{N_j^+(t)}$;
\item[-] set $V_i=-a$ for $i=k_{N_j^+(t)+1},\ldots, k_{N^c_j(t)}$;
\end{itemize}
\item otherwise set $V_i=V_i^0$.
\end{enumerate}
\end{enumerate}
\end{algorithm}

\begin{remark}
As a side result of our derivation, we constructed a Monte Carlo method for a general scalar conservation law of the form \eqref{eqch3:cls}. As for the construction of deterministic relaxation schemes for conservation laws~\cite{JinXin}, this is easily obtained by taking the limit case $\e\to 0$ in the schemes just described (see~\cite{Seaid}). 
\end{remark}


\begin{figure}[!tbp]
\centering
\includegraphics[width=0.75\textwidth]{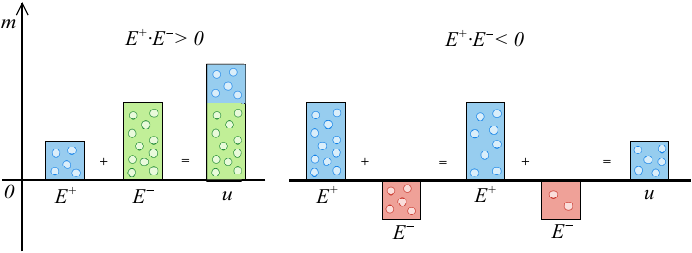}
\caption{Left: The state variable $u$ uniquely defines the number of particles in the two equilibrium states $E^{\pm}$ when the latter have the same sign. Right: There exist infinite possible ways to associate particles to the equilibrium states $E^{\pm}$ when equilibria have discordant sign.}
\label{fig.sketch_mass}
\end{figure}
\subsubsection{The case of negative solutions}
\label{sec.MCneg}
If we  apply the method to initial conditions and solutions that might be also negative, we need to associate to each particle a possible negative weight $m_i \in \{-m,m\}$, $m>0$ that we will still refer to as the particle mass, and reconstruct the solution as
\be
u_j(t)=\frac{1}{N}\sum_{k=1}^N m_k \Phi_{\Delta x}(x_{j+1/2}-X_k),\qquad
j=\ldots,-2,-1,0,1,2,\ldots
\label{eq:rhor_m}
\ee

In this context, we observe that the equilibrium states $E^{\pm}(u)$ could result either positive or negative and, consequently, also the probabilities $p^+$ and $p^-$ defined in \eqref{ch4:nd} may become negative so that the probabilistic interpretation fails. To avoid this, we can rely on the probabilistic interpretation associated to the notion of area defined by the absolute value of the equilibrium state in the computational cell.
More precisely, in each cell, we define the equilibrium probabilities as 
\be
{p^+_E}(x,\Delta t) = \frac{|E^{+}(u)|}{|{E^{+}(u)}|+|{E^{-}(u)}|} ,\quad
{p^-_E}(x,\Delta t) = \frac{|E^{-}(u)|}{|{E^{+}(u)}|+|{E^{-}(u)}|}. 
\label{pE_abs}
\ee
Note that when applying this definition,
\[\mathrm{if} \quad {E^{+}} {E^{-}} > 0 \quad \Rightarrow \quad |{E^{+}(u)}|+|{E^{-}(u)}|=|{u}|,\]
and so the probabilistic interpretation of Algorithm \ref{alg.MCscalar} remains valid. However, an indeterminate problem arises if the equilibrium states have different signs
\[\mathrm{if} \quad {E^{+}} {E^{-}} < 0 \quad \Rightarrow \quad |{E^{+}(u)}|+|{E^{-}(u)}|>|{u}|,\] 
because there is an infinite number of possible allocations of particles with opposite masses (see Figure \ref{fig.sketch_mass} for a sketch). 

To overcome this problem, after the relaxation step based on the probabilities \eqref{pE_abs}, we can follow two strategies in each cell: 
\begin{enumerate}
\item[$(i)$] keep the particle mass fixed and introduce a variable number of particles;
\item[$(ii)$] keep the number of particles fixed and introduce a variable mass. 
\end{enumerate}
Both strategies have pros and cons and we will briefly describe the different implementations below. The first choice $(i)$ would result in the need to discard or re-sample particles in each cell, and this can be done by a particle killing and replication strategy, as in the case of the source term discussed at the end of Section \ref{sec:GRW}. Precisely, in each cell $j$ of width $\Delta x$, we define the new particle numbers 
\be
{\tilde N}^+_j = {\rm Sround}\left(\frac{|E^+(u_j)|N\Delta x}{m N_j^+}\right),\qquad {\tilde N}^-_j = {\rm Sround}\left(\frac{|E^-(u_j)|N\Delta x}{m N_j^-}\right),
\ee
where $N_j^+$ and $N_j^-$ are the number of particles with positive and negative velocities in the cell, respectively. Then, we resample ${\tilde N}^+_j - N_j^+$ particles with positive velocity if ${\tilde N}^+_j > N_j^+$, or discard $N_j^+-{\tilde N}^+_j$ particles with positive velocity if ${\tilde N}^+_j < N_j^+$. The same, of course, is done for the particles with negative velocity.

In the second alternative $(ii)$ we fix the number of particles $N_j$ belonging to the $j$-th cell and, consequently, reassign their mass, which will thus depend on the cell. 
This can be done, by defining an updated mass value $\tilde{m}_j$ in each cell in order to respect the following balance:
\be
\tilde{m}_j(t) N_j(t) = \left(|{E^{+}(u_j)}|+|{E^{-}(u_j)}|\right)N \Delta x.
\label{mass_reassign_MC}
\ee
Finally, in both cases $(i)$ and $(ii)$, the sign of mass is assigned to each particle according to the sign of the corresponding equilibrium. The details of the second strategy $(ii)$ are reported in Algorithm \ref{alg.MCneg}, where for simplicity we restrict to the limiting case $\e \to 0$.

\begin{algorithm}[h!]
\label{alg.MCneg}
\caption{Weighted Monte Carlo for scalar conservation laws}
\begin{enumerate}
\item Compute
$X_i=X^{0}_i+V^{0}_i \Delta t, \quad i=1,\ldots,N$.
\item For each grid cell $j$ reconstruct $u_j(t)$ from \eqref{eq:rhor_m} and compute $\tilde{m}_j$ using \eqref{mass_reassign_MC}.
\begin{enumerate}
\item if ${E_j^{+}(u_j)}>0$, set $m_j^+ = \tilde{m}_j$ otherwise set $m_j^+ = -\tilde{m}_j$;
\item if ${E_j^{-}(u_j)}>0$, set $m_j^- = \tilde{m}_j$ otherwise set $m_j^- = -\tilde{m}_j$;
\end{enumerate}
\item For each sample $(X_i,V^{0}_i)$ in cell $j$,
\begin{enumerate}
\item with probability ${p^+_E}(x_j,t)$ given by \eqref{pE_abs} set
$V_i=a$, $m_i=m^+_j$,
\item otherwise set $V_i=-a$, $m_i=m^-_j$.
\end{enumerate}
\end{enumerate}
\end{algorithm}

\subsection{The Gradient-based Monte Carlo method}
In this section we show how the Monte Carlo approach just described can be improved significantly using a gradient random walk strategy. In fact, by introducing the auxiliary variables $w=\partial u/\partial x$ and $z=\partial v/\partial x$, we can rewrite (\ref{eqch3:jx}) in the form
\be
\begin{aligned}
\frac{\partial w}{\partial t} + \frac{\partial z}{\partial x} &=  0,\\
\frac{\partial z}{\partial t} + a^2 \frac{\partial w}{\partial x} &= -\frac1{\e} \left(z-F'(u)w\right).
\end{aligned}
\label{eqch3:jxg}
\ee

Most of the arguments presented for system (\ref{eqch3:jx}) apply again. First we introduce the new variables
\[
g^+=\frac{(aw+z)}{2a},\qquad g^-=\frac{(aw-z)}{2a},
\]
which satisfy the diagonal system
\be
\begin{aligned}
\frac{\partial g^+}{\partial t} + {a}\frac{\partial g^+}{\partial x} &= -\frac1{\e} (g^+-D^+(u,w)),\\
\frac{\partial g^-}{\partial t} - {a} \frac{\partial g^-}{\partial x} &= -\frac1{\e} (g^--D^-(u,w)),
\end{aligned}\label{ch4:gt2}
\ee
where the equilibrium states now read
\[
D^+(u,w)=\frac{w({a}+F'(u))}{2{a}}, \qquad D^-(u,w)=\frac{w({a}-F'(u))}{2{a}}.
\]
After splitting the above system, the transport step and the collision/relaxation step can be solved using the same Monte Carlo approach described in the previous sections, but with a substantial difference: there is no need to reconstruct the solution $u(x,t)$ on a space grid during the relaxation step. In fact, given a set of samples located in $X_1$,$\ldots$, $X_N$ with positive and negative masses $m_i \in \{-m,m\}$ we can compute $u(X_i)$ following \eqref{eqch3:cs3}. 
Since  $D^+(u,w)+D^-(u,w)=w$, under the subcharacteristic condition $a>|F'(u)|$, the probabilities of a random velocity change read
\be
p^+_D(x,\Delta t)=\frac{D^{+}(u,w)}{w}=\frac{a+ F'(u)}{2a}, \qquad 
p^-_D(x,\Delta t)=\frac{D^{-}(u,w)}{w}=\frac{a- F'(u)}{2a}.
\ee 
Therefore, applying the gradient random walk strategy, the probabilities turn out to be dependent only on $u$ (and not $w$), which is now a particle-dependent (and not grid-dependent) variable. This aspect is of utmost importance because leads to deal with a grid free method and consequent advantages both in terms of accuracy and computational efficiency.

Thus, we can solve system \eqref{eqch3:jxg} with a Gradient-based Monte Carlo method (GBMC) as follows. Starting with a set of samples $(X^0_1,V^0_1), \ldots, (X^0_{N},V^0_N)$, with $V^0_i\in\{-a,a\}$ and each sample has mass $m_i\in\{-m,m\}$, a new set of samples $(X_1,V_1),\ldots,(X_{N},V_N)$ is generated as reported in Algorithm \ref{alg.RWMCscalar}. 

\begin{algorithm}[ht]
\label{alg.RWMCscalar}
\caption{Gradient-based Monte Carlo for $2\times 2$ hyperbolic relaxation systems}
\begin{enumerate}
\item Compute
$X_i=X_i^0+V_i^0 \Delta t,\quad i=1,\ldots,N$.
\item For each sample $X_i$, evaluate $u_i(t)$ using \eqref{eqch3:cs3}.
\item For each sample pair $(X_i,V_i^0)$,
\begin{enumerate}
\item with probability $1-e^{-\Delta t/\e}$ do the following
\begin{itemize}
\item[-] with probability $\displaystyle\frac{a+F'(u_i(t))}{2a}$ set
$V_i=a$, 
\item[-] otherwise set $V_i=-a$.
\end{itemize}
\item otherwise set $V_i=V_i^0$.
\end{enumerate}
\end{enumerate}
\end{algorithm}

Clearly, similar to the standard Monte Carlo case, taking the limit $\e\to 0$ in the above algorithm leads to a Gradient-based Monte Carlo method for a general scalar conservation law of the form \eqref{eqch3:cls}. Let us mention that other approximations for scalar conservation laws based on Gradient-based Monte Carlo are found in the literature for rather specific situations, like the Burgers equation (see~\cite{Rob} for details).

\begin{remark}
\label{rmk.REC}

In addition to left reconstruction \eqref{eqch3:cs3} denoted by $u_i^L(t)$ or the analogous right reconstruction $u_i^R(t)$, other reconstructions can be implemented that interpolate between the two. For example, to avoid asymmetric reconstructions, a weighted average of left and right reconstructions can be computed for the whole computational range $[x_{\min},x_{\max}]$, where $x_{\min}={\min_i}\{X_i\}$, $x_{\max}={\max_i}\{X_i\}$ considering linearly distributed weights in $[0,1]$:
\be
u_i(t) = (1-w_i) u_i^L(t) + w_i u_i^R(t),\qquad \omega_i = (X_i-x_{\min})/(x_{\max}-x_{\min}),
\label{eq:reco}
\ee 
Of course, more sophisticated reconstruction can be designed similarly. In the numerical examples, if not otherwise stated, we will make use of \eqref{eq:reco}.
\end{remark}

\subsection{Numerical examples for scalar conservation laws}
In this section we compare the numerical solutions obtained with the standard Monte Carlo approach and the Gradient-based Monte Carlo method for different problems governed by scalar conservation laws including an empirical convergence rate test. Indeed, since the $\e\to 0$ case is the most challenging as the limiting nonlinear scalar conservation law can form discontinuities in finite time, we will limit the presentation of results to this situation.


\begin{figure}[!tbp]
\centering
\includegraphics[width=0.42\textwidth]{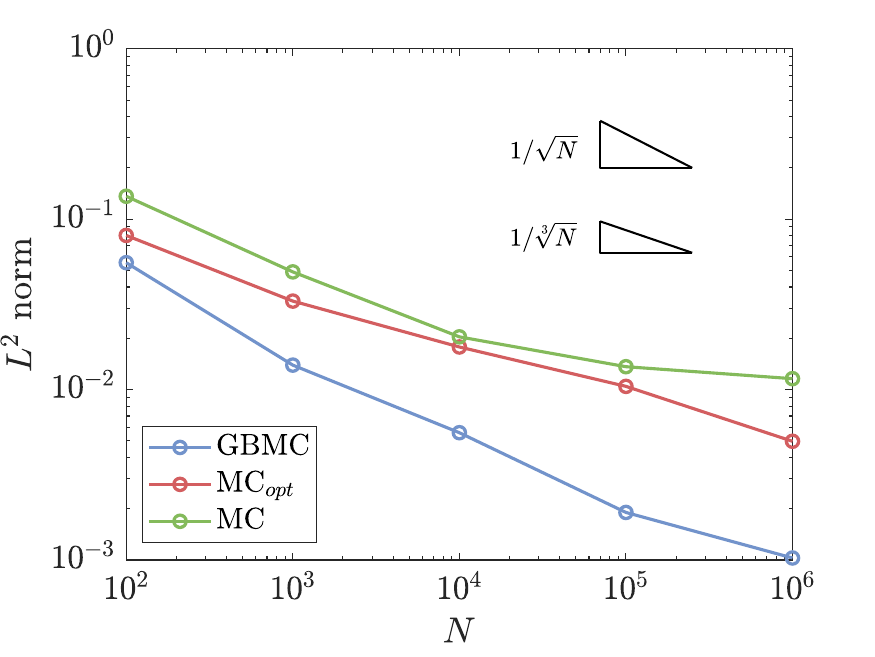}
\includegraphics[width=0.42\textwidth]{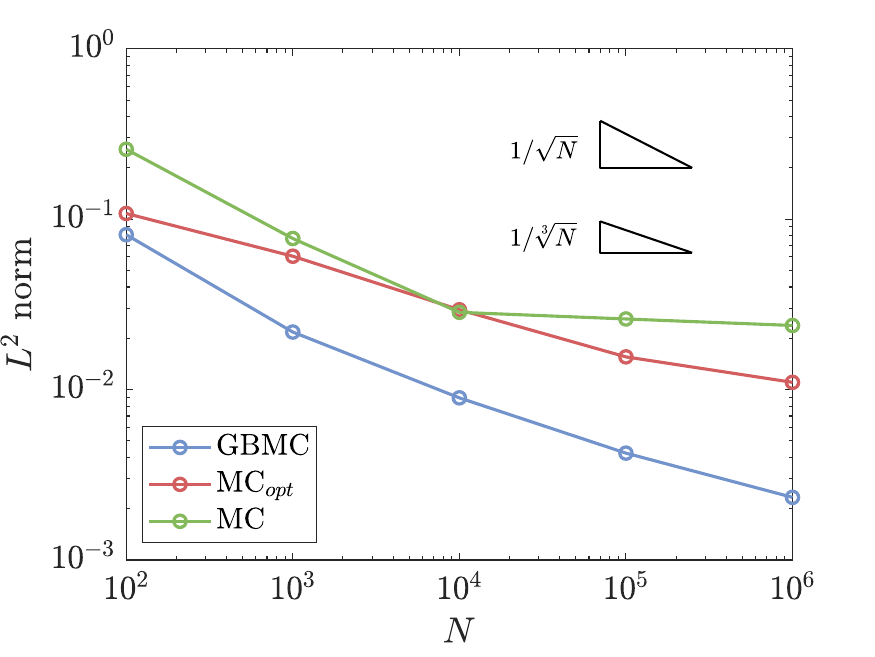}
\caption{Convergence rate analysis. Comparison of the relative $L^2$ error norms of the direct Monte Carlo (MC) and the Gradient-based Monte Carlo (GBMC) with respect to the number of particles $N$ for the solution of the inviscid Burgers equation at $t=2.5$ with normal distribution as initial datum (left) and at $t=0.5$ with sinusoidal distribution as initial datum (right).}
\label{convergence_rate}
\end{figure}
\begin{table}[b]
\label{tab.convergence}
\centering
\begin{tabular}{c | c c | c c }
\toprule
	 \multirow{2}{*}{$N$} &\multicolumn{2}{c|}{$L^2_{MC}/L^2_{GBMC}$} &\multicolumn{2}{c}{$L^2_{MC_{opt}}/L^2_{GBMC}$}\\
	 &(a) &(b) &(a) &(b)\\
\midrule
	$10^2$ &2.45 &3.17 &1.45 &1.33\\
	$10^3$ &3.52 &3.54 &2.38 &2.79\\
	$10^4$ &3.65 &3.18 &3.18 &3.29\\
	$10^5$ &7.16 &6.13 &5.49 &3.67\\
	$10^6$ &11.27 &10.20 &4.84 &4.73\\
\bottomrule
\end{tabular}
\caption{Convergence rate analysis. Ratio $L^2_{MC}/L^2_{GBMC}$ with respect to the number of particles $N$ for the solution of the inviscid Burgers equation at $t=2.5$ with normal distribution as initial datum (a) and at $t=0.5$ with sinusoidal distribution as initial datum (b).}
\end{table}

\subsubsection{Empirical convergence rate} 
In the first test case, we compute the convergence rate of the methods with respect to the number of particles solving the inviscid Burgers equation, corresponding to $F(u)=u^2/2$ in \eqref{eqch3:jx} and $\e\to 0$. First, we consider a normal distribution as initial datum and run the simulations up to $t=2.5$, prior to the classical shock formation; then, a test with sinusoidal initial condition (i.e., involving also particles with negative masses) with solution at $t=0.5$ is taken into account. The comparison of the convergence rates obtained in the two test cases with 
\begin{itemize}
\item the Monte Carlo method while keeping the mesh size $\Delta x$ fixed with the particles number refinement (MC), 
\item the Monte Carlo method with optimal choice of the mesh size $\Delta x$ in function of the $N$ particles amount (MC$_{opt}$), as discussed in Appendix \ref{appendix_MC},
\item the Gradient-based Monte Carlo method (GBMC),
\end{itemize}
are shown in Figure \ref{convergence_rate} with respect to the relative $L^2$ norm. Here, we considered as reference solution the one obtained with a finite volume Godunov method with a very refined mesh grid. Both MC, MC$_{opt}$ and GBMC errors are computed with respect to the the mean of 5 runs with $N$ particles, while keeping an empirically chosen sufficiently small time step $\Delta t$ fixed (to prevent the temporal error from overriding the error due to the choice of particles number and, thus, alter the convergence curves). As discussed in Appendix \ref{appendix}, it can be observed that, both the standard MC and GBMC methods present an error decay of $O(1/\sqrt{N})$, while when optimizing the choice of the mesh size, the MC method converges with $O(1/\sqrt[3]{N})$. Nevertheless, the error produced by the GBMC is always smaller than that obtained by applying either the optimized or non-optimized MC method. In particular, in Table \ref{tab.convergence}, we report the precise gain factor from using the GBMC instead of both the MC versions in terms of accuracy as a function of the number o particles $N$, intended as ratio of the two norms $L^2_{MC}/L^2_{GBMC}$ and $L^2_{MC_opt}/L^2_{GBMC}$. Let su finally point out that in both the test cases the decay of the order of accuracy in the MC method observed for large $N$ is due to the fact that the mesh error dominates the statistical error due to the particle number. This is in good agreement with the analysis in Appendix \ref{appendix}.

\begin{figure}[!htbp]
\centering
\includegraphics[width=0.35\textwidth]{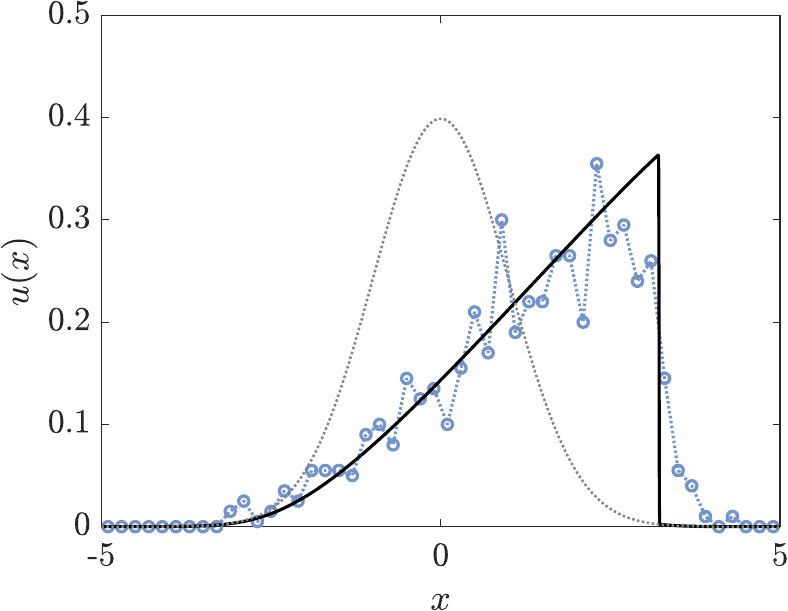}
\includegraphics[width=0.35\textwidth]{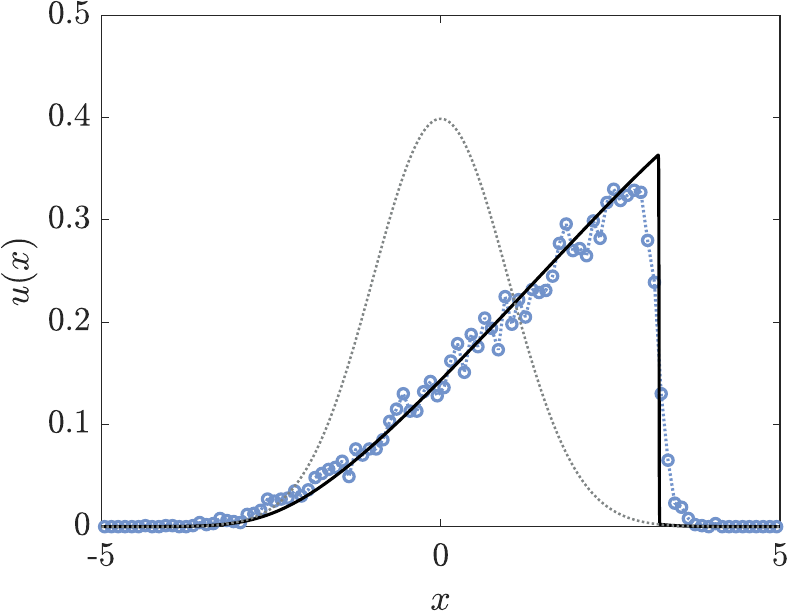}
\includegraphics[width=0.35\textwidth]{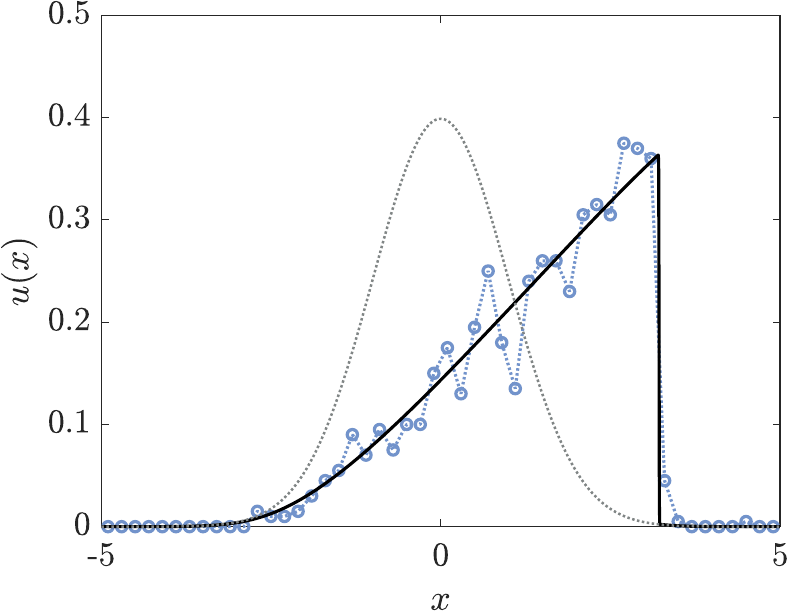}
\includegraphics[width=0.35\textwidth]{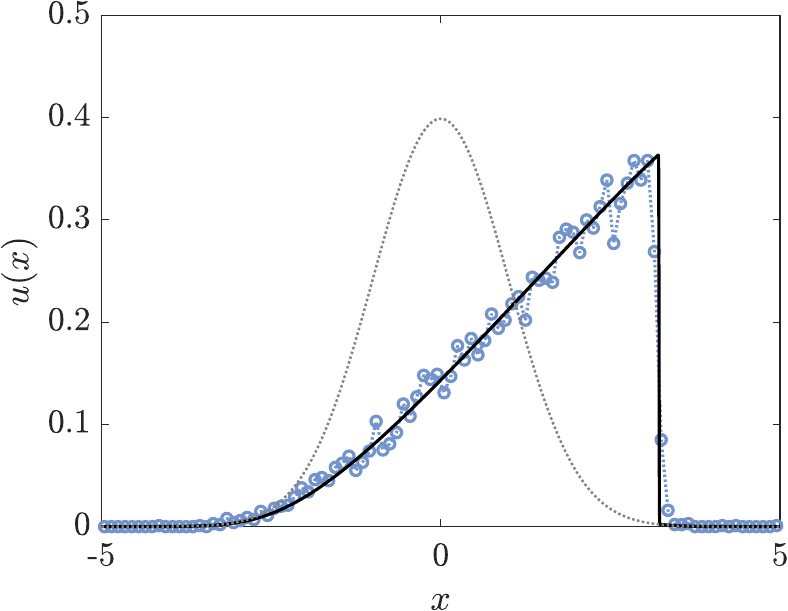}
\includegraphics[width=0.35\textwidth]{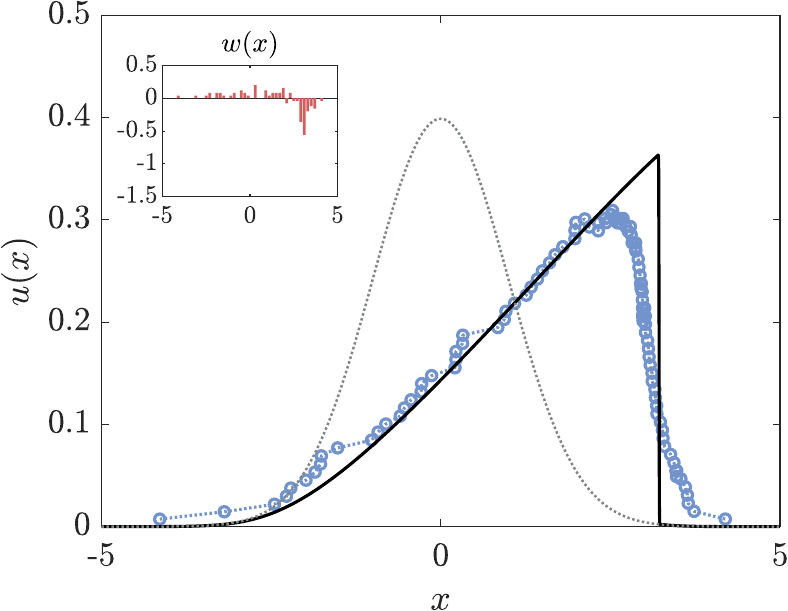}
\includegraphics[width=0.35\textwidth]{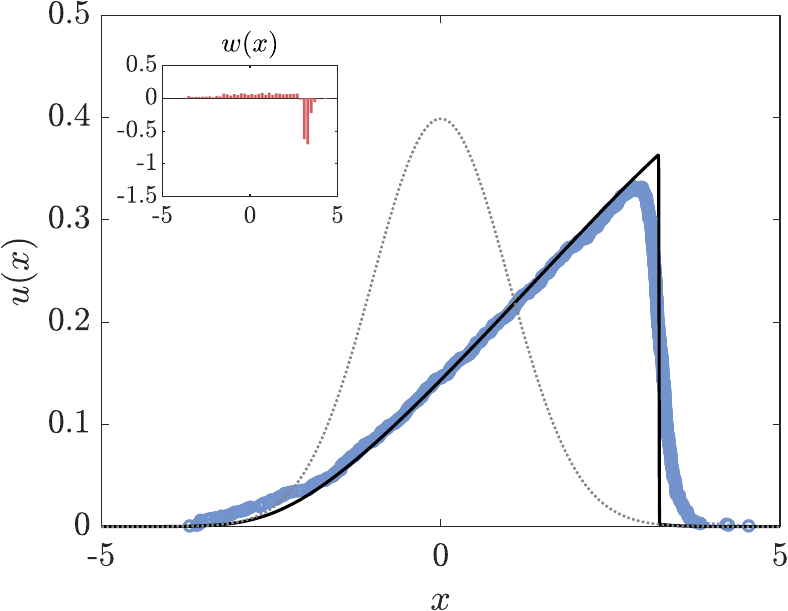}
\includegraphics[width=0.35\textwidth]{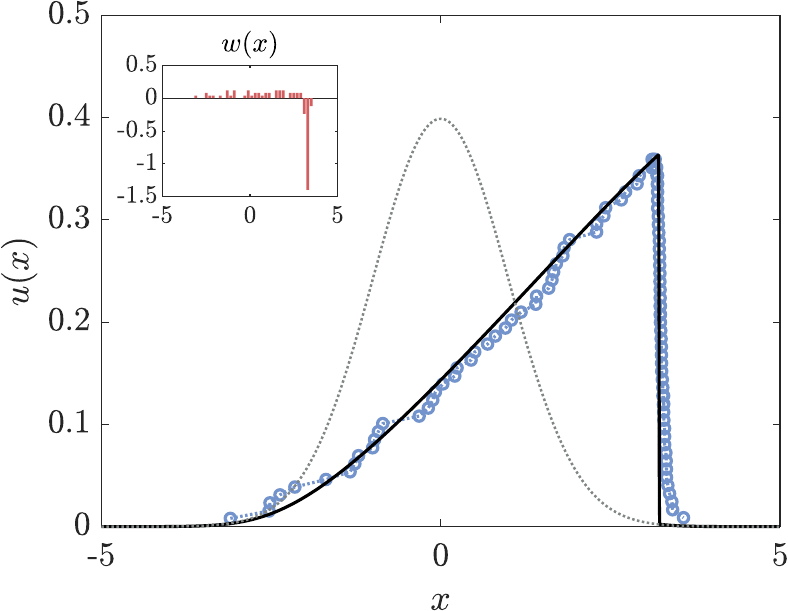}
\includegraphics[width=0.35\textwidth]{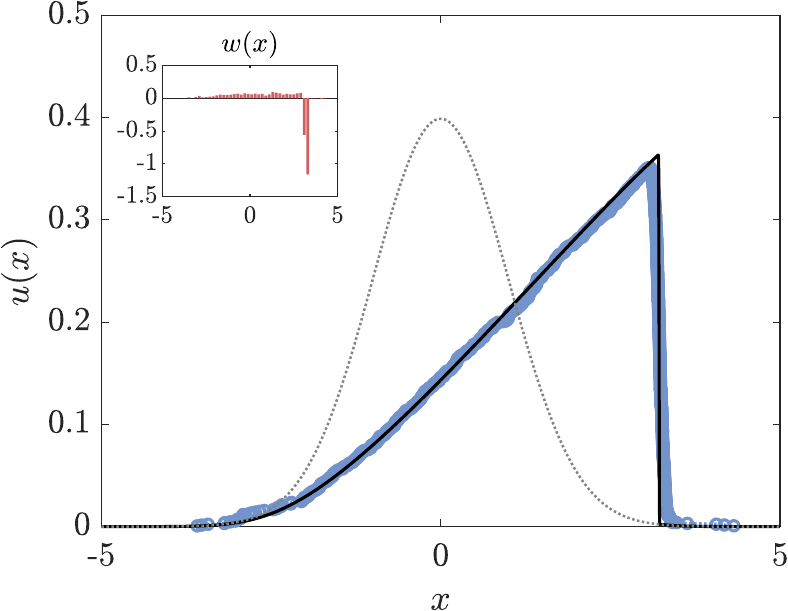}
\caption{Test 1(a), Burgers equation. Solution at $t=10$ with a Gaussian density as initial datum, shown in gray dotted line. Here $a=0.4$ and $\Delta t = 0.5$ (first and third row, left), $\Delta t = 0.25$ (first and third row, right), or $\Delta t = 0.1$ (second and fourth row). First two rows: Standard MC with $N=1000$ particles and $M=50$ grid points (left), and $N=10000$ and $M=100$ grid points (right). Last two rows: GBMC method with $N=100$ (left) and $N=1000$ (right). The subplots refer to the histogram of the space derivative of the state variable. The reference solution is reported in black solid line.}
\label{burgers_gaussIC}
\end{figure}

\begin{figure}[!tb]
\centering
\includegraphics[width=0.42\textwidth]{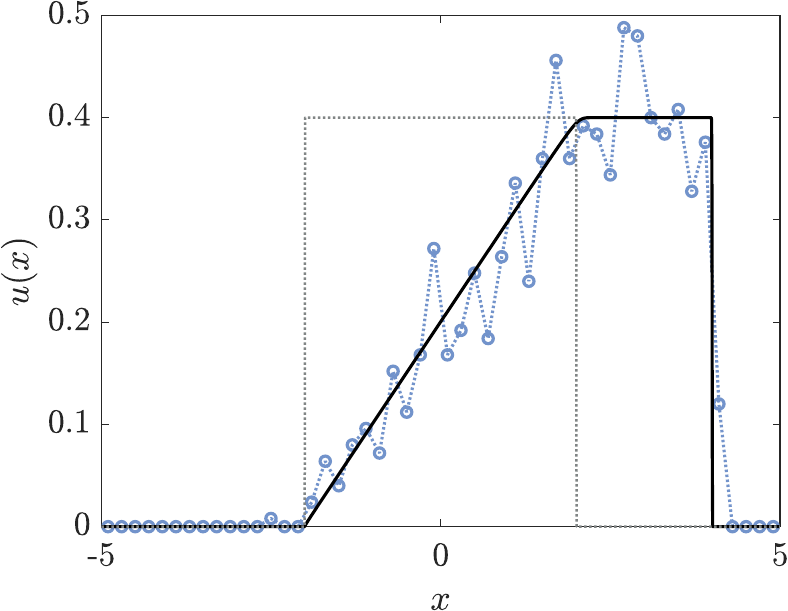}
\includegraphics[width=0.42\textwidth]{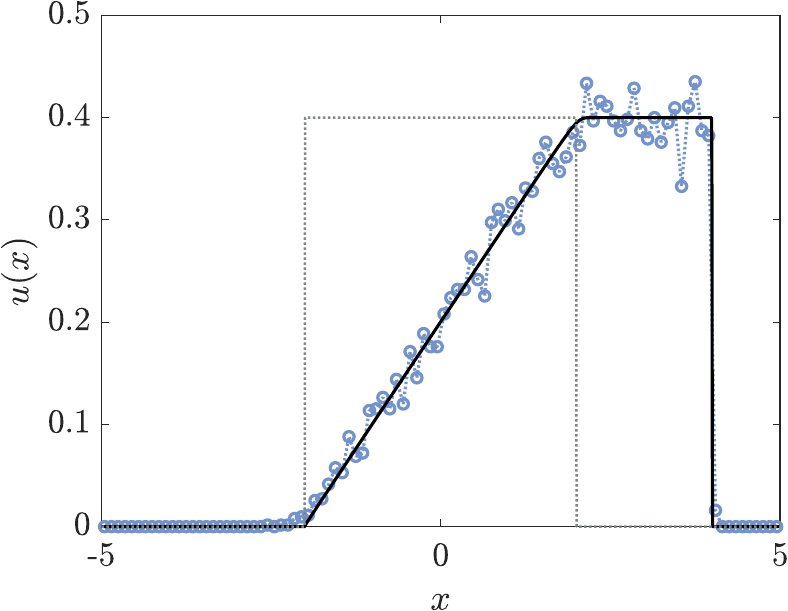}
\includegraphics[width=0.42\textwidth]{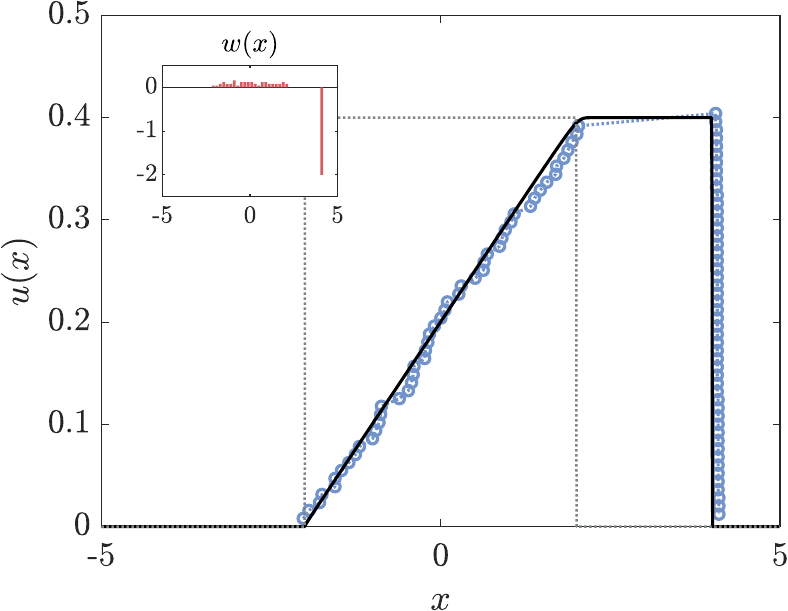}
\includegraphics[width=0.42\textwidth]{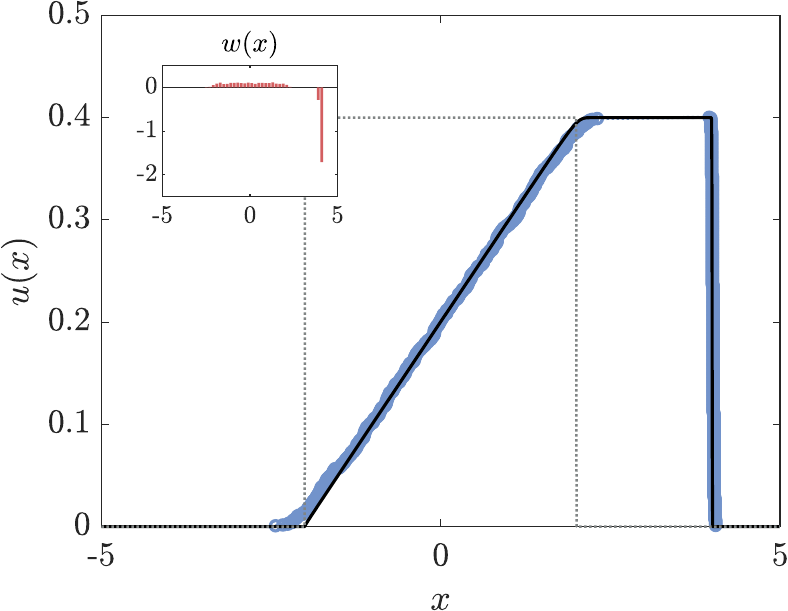}
\caption{Test 1(b), Burgers equation. Solution at $t=10$ with a square wave as initial datum, shown in gray dotted line. Here $a=0.6$ and $\Delta t = 0.01$. First row: Standard MC with $N=1000$ particles and $M=50$ grid points (left), and $N=10000$ and $M=100$ grid points (right). Last row: GBMC method with $N=100$ (left) and $N=1000$ (right). The subplots refer to the histogram of the space derivative of the state variable. The reference solution is reported in black solid line.}
\label{burgers_squareIC}
\end{figure}

\begin{figure}[!tb]
\centering
\includegraphics[width=0.42\textwidth]{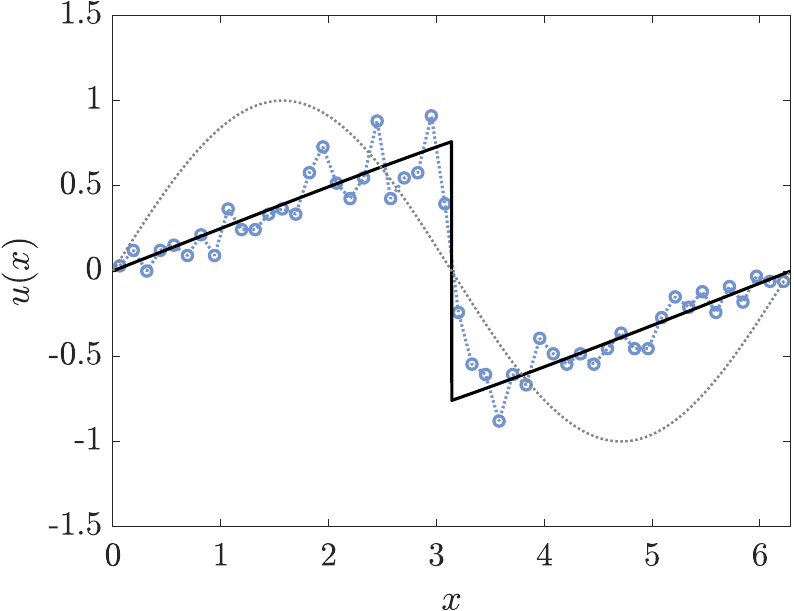}
\includegraphics[width=0.42\textwidth]{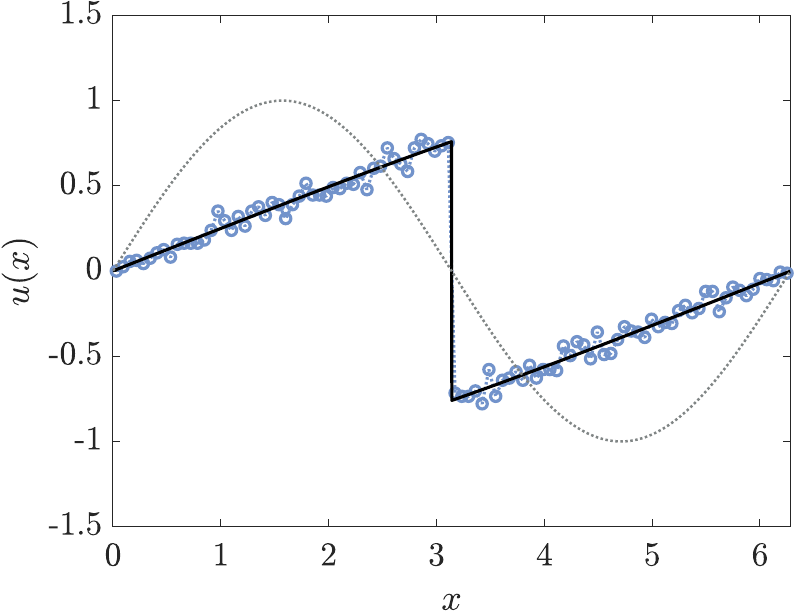}
\includegraphics[width=0.42\textwidth]{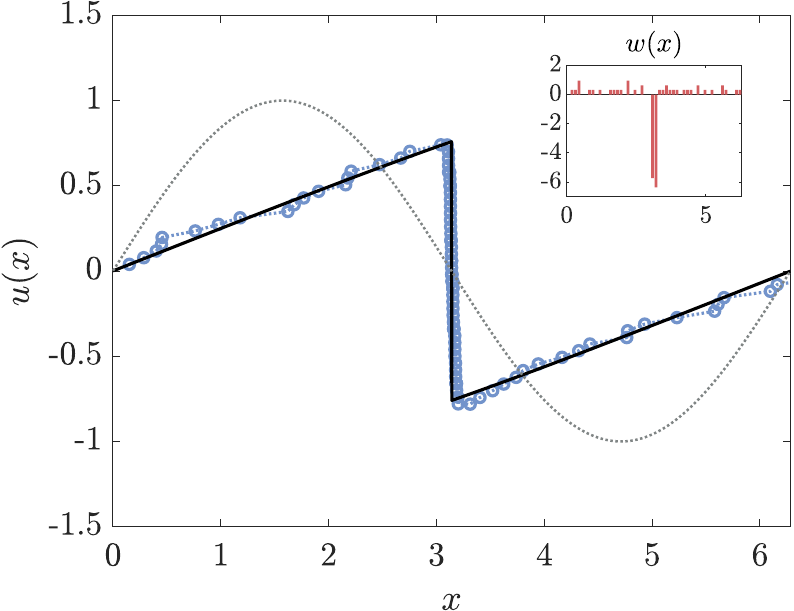}
\includegraphics[width=0.42\textwidth]{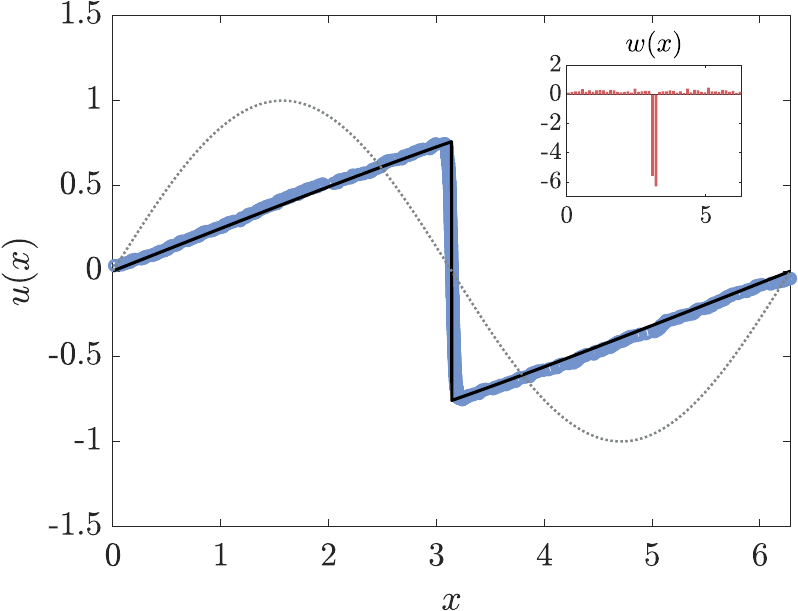}
\caption{Test 1(c), Burgers equation. Solution at $t=3$ with a sinusoidal initial datum, shown with gray dotted line. 
Here $a=1.5$ and $\Delta t = 0.01$. First row: Standard MC with $N=1000$ particles and $M=50$ grid points (left), and $N=10000$ and $M=100$ grid points (right). Last row: GBMC method with $N=100$ (left) and $N=1000$ (right). The subplots refer to the histogram of the space derivative of the state variable. The reference solution is reported in black solid line.}
\label{burgers_sinIC}
\end{figure}


\subsubsection{Inviscid Burgers equation} 
The MC (without low variance technique) and GBMC methods are then applied again to the inviscid Burgers equation considering four different initial conditions.

\begin{itemize}
\item\textbf{Test 1(a):} In the first case, the initial datum of the inviscid Burgers equation is a Gaussian density with zero mean and unit variance.

\item\textbf{Test 1(b):} In the second case, we consider a rectangular wave as initial datum  
\be
u(x,0)= \left\{
\begin{array}{ll}
 0.4 \quad &{\rm if}\,\, -2\leq x\leq2\\ 
 0 \quad &{\rm otherwise}\,.
 \end{array}\right.
\ee

\item\textbf{Test 1(c):} In the third case, we fix as initial condition a sinusoidal function 
\be
u(x,0)=\sin(x)\,,
\ee
to assess the performance of the methods even when considering negative solutions, hence introducing particles with negative mass.

\end{itemize}

For each test case, results are given at equal time discretization for the two methods in Figures~\ref{burgers_gaussIC}--\ref{burgers_sinIC}. With the Monte Carlo method, $N=1000$ or $N=10000$ particles and $M=100$ cells for the domain discretization are considered; while with the GBMC approach, $N=100$ or $N=1000$ particles. The reference solution is obtained employing a finite volume Godunov method with a very refined spatial grid.
By comparing results obtained with the two methods, the remarkable improvement in variance reduction and in capturing the shock fronts obtained by the usage of the GBMC method appears evident, even if considering a largely reduced amount of particles with respect to the standard MC. Moreover, in Figure~\ref{burgers_gaussIC} it can also be observed the influence of the choice of the time step size on the final solution for both methods. In particular, we note that the choice of $\rm{CFL} =1$ (and thus to $\Delta t = 0.5$ for a mesh with $M=50$ and $\Delta t = 0.25$ for a grid with $M=100$ cells), leads to an evident numerical dissipation, visible near the shock. In fact, we remark that the Monte Carlo approach here proposed, for $N \to \infty$ and $\rm{CFL} = 1 \Rightarrow \Delta t = \Delta x/a$ coincides with the Lax--Friedrichs scheme (see~\cite{JinXin} for further details), which is known to produce high numerical dissipation when $\rm{CFL} = 1$.

\begin{figure}[t!]
\centering
\includegraphics[width=0.42\textwidth]{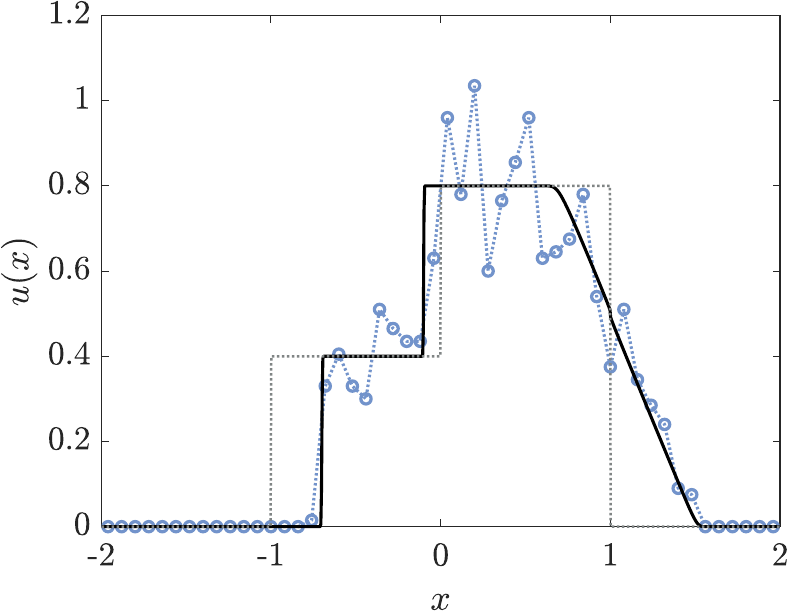}
\includegraphics[width=0.42\textwidth]{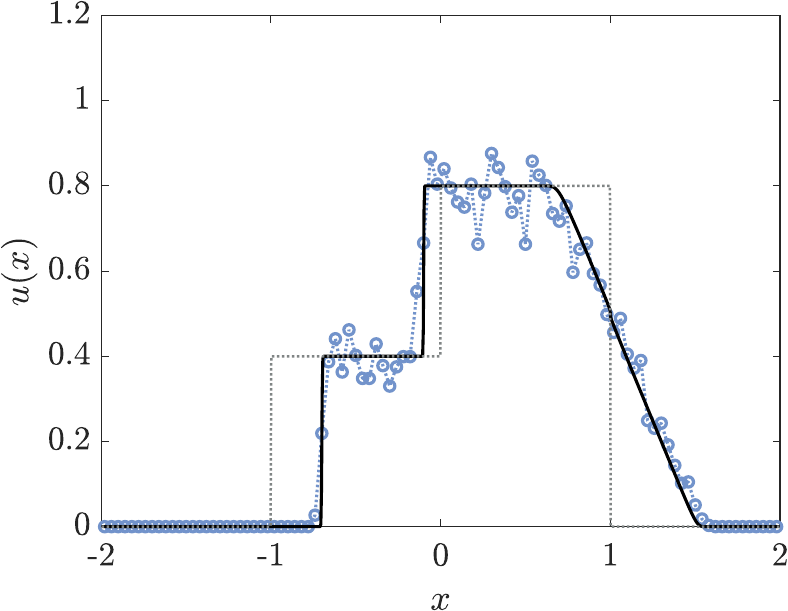}\\
\includegraphics[width=0.42\textwidth]{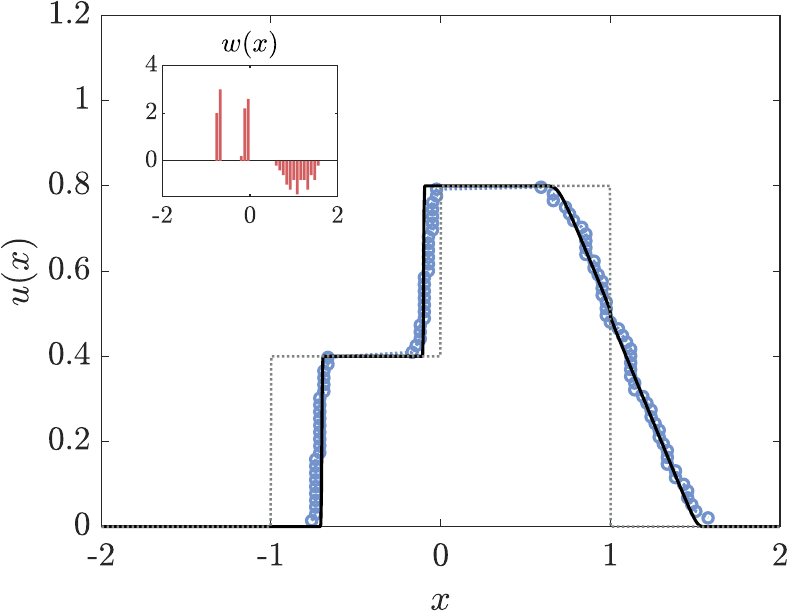}
\includegraphics[width=0.42\textwidth]{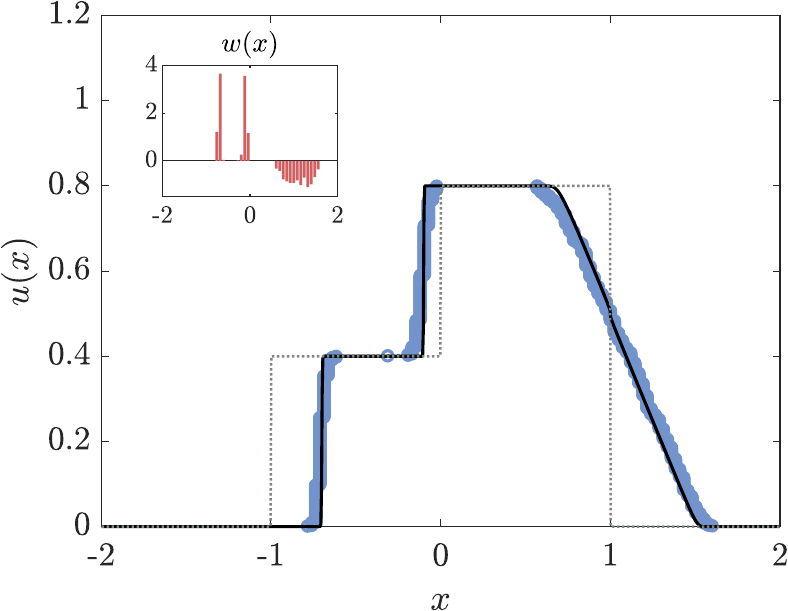}
\caption{Test 2, LWR traffic model. Solution at $t=0.5$ of the Riemann problem \eqref{eq:RP}. 
Initial datum shown as gray dotted lines. 
Here $a=1.2$, $\Delta t = 0.01$. First row: MC with $N=1000$ particles and $M=50$ grid points (left), and $N=10000$ and $M=100$ (right). Last row: GBMC  with $N=100$ (left) and $N=1000$ (right). The subplots refer to the histogram of the space derivative of the state variable. The reference solution is reported in black solid line.}
\label{LWR_RP}
\end{figure}

\subsubsection{Lighthill-Whitham-Richards traffic model} 
As a second scalar conservation law, we consider the Lighthill-Whitham-Richards (LWR) traffic model, for which $F(u)=u-u^2$ in \eqref{eqch3:jx} and $\e\to 0$, taking into account the following initial conditions.
\begin{itemize}
\item
\textbf{Test 2:} We consider the Riemann Problem (RP) presented in~\cite{difrancesco2017}, in which
\be
u(x,0) = \left\{
                \begin{array}{ll}
                0.4 \quad &\mathrm{if} \,\,-1 \le x \le 0,\\
                0.8 \quad &\mathrm{if} \,\,0 < x \le 1,\\
                0 \quad &\mathrm{otherwise}.
                \end{array}
            \right.
            \label{eq:RP}
\ee

\end{itemize}

Similar observations to those for Test 1 can also be made here when comparing the results in Figure \ref{LWR_RP}, which shows numerical solutions of the LWR test cases obtained applying the standard Monte Carlo and the GBMC method considering different amounts of particles. The well reduced variance of GBMC solutions can again be appreciated even if fewer particles are used than in direct Monte Carlo, which also brings a consequent speeding up of the simulation. Moreover, it can be observed the capability of the proposed method to well capture the sharp discontinuities arising in the dynamics, thanks to its adaptive nature of following the solution through particles especially where the gradient is large.
\section{Extensions to hyperbolic systems of conservation laws}
\label{sec:GBMCsy}
In this section we show how to extend the Monte Carlo and Gradient-based Monte Carlo techniques to the case of hyperbolic relaxation approximations to systems of conservation laws. As before, the methods we derive here work uniformly with respect to the stiff relaxation rate and, in the zero relaxation limit, originate a Monte Carlo or Gradient-based Monte Carlo method for the corresponding hyperbolic system of conservation laws. Our attention will be focused in particular on the behavior of methods in that limit.  

\subsection{A Monte Carlo approach}
Consider the system of conservation laws in one space variable
\be
\frac{\partial \U}{\partial t} + \frac{\partial \FF(\U)}{\partial x} = 0,
\label{eq:scons} 
\ee
with $x\in\Omega\subseteq \mathbb{R}$, $\U=(u_1,\ldots,u_n)\in\mathbb{R}^n$, $n\geq 2$. The above system is strictly hyperbolic if the Jacobian matrix  $\FF'(\U)$ admits $n$ distinct real eigenvalues $\lambda_1 <\ldots < \lambda_n$. The system is complemented with the initial conditions $\U(x,0)=\U_0(x)$ and suitable boundary conditions. 


The relaxation approximation now reads~\cite{JinXin}
\be
\begin{aligned}
\frac{\partial \U}{\partial t} + \frac{\partial \V}{\partial x} &=  0,\\
\frac{\partial \V}{\partial t} + A^2 \frac{\partial \U}{\partial x} &= -\frac1{\e} (\V-\FF(\U)),
\end{aligned}
\label{eq:reln}
\ee
where $\V\in\mathbb{R}^n$ and $A^2={\rm diag}\{a_1^2,\ldots,a_n^2\}$ must satisfy the dissipative condition $A^2 > \FF'(\U)^2$ for all $\U$, with initial conditions $\V(x,0)=\V_0(x)$ and defined boundary conditions. Notice that, for $\U$ varying in a bounded domain, the dissipative condition can always be satisfied by choosing sufficiently large $A$, but because of the CFL condition, for numerical stability, it is desirable to obtain the smallest $A$ meeting the criterion. Following~\cite{JinXin}, we say that system \eqref{eq:reln} is dissipative if it is strictly stable in the sense of Majda-Pego, which is satisfied if $a^2_h > \lambda^2_h,\, h=1,\ldots,n$.

The diagonal variables are
\[
\F^{\pm}=A^{-1}\frac{A\U\pm \V}{2},
\]
which yield the system
\be
\begin{aligned}
\frac{\partial \F^+}{\partial t} + A\frac{\partial \F^+}{\partial x} &=  -\frac1{\e} (\F^+-\E^+(\U))\\
\frac{\partial \F^-}{\partial t} - A \frac{\partial \F^-}{\partial x} &= -\frac1{\e} (\F^--\E^-(\U)),
\end{aligned}
\label{eq:relnd}
\ee
with
\be
\E^{\pm}(\U) = A^{-1}\frac{A\U\pm \FF(\U)}{2}. 
\ee
As in the case of scalar of conservation laws, if we assume that initial conditions might be also negative, the equilibrium states $\E^{\pm}(\U)$ could result either positive or negative as well. Therefore, recalling the discussion presented in Section \ref{sec.MCneg}, we define the probabilities of a random velocity change for each family of particles $h$ as 
\be
p_E^{h,+}(x,\Delta t) = \frac{|E_h^{+}(\U)|}{|{E_h^{+}(\U)}|+|{E_h^{-}(\U)}|}, \qquad
p_E^{h,-}(x,\Delta t) = \frac{|E_h^{-}(\U)|}{|{E_h^{+}(\U)}|+|{E_h^{-}(\U)}|},
\label{pE_syst}
\ee
and proceed similarly to the scalar case by either considering a variable particle number or  including the update of the particles' mass. In the latter case, in each cell $j$ of width $\Delta x$ we define for each species $h$ the mass value that fulfills the following relation:
\be
\tilde{m}^h_j N_j^h = \left(|{E_h^{+}(\U_j)}|+|{E_h^{-}(\U_j)}|\right) {\Delta x \,N}.
\label{mass_reassign_MCsyst}
\ee
Here $N_j^h$ is the number of particles of the family $h$ inside the $j$-th cell. In Algorithm \ref{alg.MCsyst} we summarize the weighted Monte Carlo scheme in the limiting case $\e\to 0$
starting with $n$ sets of samples
$(X^{h,0}_1,V^{h,0}_1)$, $\ldots$, $(X^{h,0}_{N^h},V^{h,0}_{N^h})$, with mass $(m^{h,0}_1,\ldots,,m^{h,0}_{N^h})$, $h=1,\ldots,n$, where
$V^{h,0}_i\in\{-a_h,a_h\}$.

\begin{algorithm}[h!]
\label{alg.MCsyst}
\caption{Weighted Monte Carlo for systems of conservation laws}
\begin{enumerate}
\item Compute
$X^h_i=X^{h,0}_i+V^{h,0}_i \Delta t,\quad h=1,\ldots,n, \quad i=1,\ldots,N^h$.
\item For each grid cell $j$ reconstruct $u_j^h(t)$
from \eqref{eq:rhor_m} and compute $\tilde{m}^h_j$ using \eqref{mass_reassign_MCsyst}.
\begin{enumerate}
\item if ${E_h^{+}(\U_j)}>0$, set $m_j^{h,+} = \tilde{m}^h_j$ otherwise set $m_j^{h,+} = -\tilde{m}^h_j$;
\item if ${E_h^{-}(\U_j)}>0$, set $m_j^{h,-} = \tilde{m}^h_j$ otherwise set $m_j^{h,-} = -\tilde{m}^h_j$;
\end{enumerate}
\item For each sample $(X^h_i,V^{h,0}_i)$ in cell $j$
\begin{enumerate}
\item with probability ${p_E^{h,+}}(x_j,t)$ as given in \eqref{pE_syst} set
$V^h_i=a_h$, $m^h_i=m^{h,+}_j$,
\item otherwise set $V^h_i=-a_h$, $m^h_i=m^{h,-}_j$.
\end{enumerate}
\end{enumerate}
\end{algorithm}

\begin{remark}
We remark that the low variance technique presented in Algorithm \ref{alg.VarReduction} can be straightforwardly implemented also in this context. Clearly, taking the limit $\e\to 0$ in Algorithm \ref{alg.MCsyst} yields a Monte Carlo method for general systems of conservation laws. It is worth to notice that the same strategy can be adopted starting from other relaxation approximations such as the one proposed in \cite{ADN}.
\end{remark}

\subsection{The Gradient-based Monte Carlo method}
\label{sect.RW_systems}
To extend the Gradient-based Monte Carlo method to systems of conservation laws, we need to introduce the quantities $\W=\partial \U / \partial x$, $\Z=\partial \V / \partial x$ to get from (\ref{eq:reln})
\be
\begin{aligned}
\displaystyle\frac{\partial \W}{\partial t} + \frac{\partial \Z}{\partial x} &=  0,\\
\displaystyle\frac{\partial \Z}{\partial t} + A^2 \frac{\partial \W}{\partial x} &= -\frac1{\e} (\Z-\FF'(\U)\W).
\end{aligned}
\label{eq:reln2}
\ee
In diagonal form, using the variables
\[
\G^{\pm}=A^{-1}\frac{A\W\pm \Z}{2},
\]
we obtain
\be
\begin{aligned}
\frac{\partial \G^+}{\partial t} + A\frac{\partial \G^+}{\partial x} &=  -\frac1{\e} (\G^+-\D^+(\U,\W))\\
\frac{\partial \G^-}{\partial t} - A \frac{\partial \G^-}{\partial x} &= -\frac1{\e} (\G^--\D^-(\U,\W)),
\end{aligned}
\label{eq:relnd2}
\ee
where now 
\be
\D^{\pm}(\U,\W)=A^{-1}\frac{(A\pm \FF'(\U))\W}{2}.
\label{D_RWMC}
\ee
We can now observe that, in the case of general hyperbolic systems of conservation laws, we cannot make the probabilities of velocity switches in the relaxation process independent on the vector $\W$, unless we can diagonalize the Jacobian matrix $\FF'(\U)$. 
Thus, in the general case of systems it is not possible to avoid the introduction of a spatial grid. We emphasize, however, that a major difference persists from the Monte Carlo method that ensures that better accuracy is achieved: equilibrium states are defined for each individual particle, $\W$ being reconstructed in the grid cells but $\U$ being particle-dependent.  This allows the solution $\U$ to be reconstructed by cumulative distribution in a manner analogous to the scalar case.

For the sake of simplicity, let us assume to be in the situation in which the original system \eqref{eq:scons} can be written in characteristic form through the Riemann invariants, so that the matrix $\FF'(\U)$ is the diagonal matrix of the eigenvalues $\lambda_1,\ldots,\lambda_n$ and we can write equilibria \eqref{D_RWMC} component-wise for each species $h$ as 
\be
D_h^{{\pm}}(\U,\W)=\frac{(a_h\pm \lambda_h(\U))w_h}{2a_h},\qquad h=1,\ldots,n.
\ee
Under the condition $a_h > |\lambda_h(\U)|$, this permits to remove the grid dependence from the probability of random velocity changes, which read
\be
\begin{split}
{p_D^{h,+}}(x,\Delta t) &= \frac{D_h^{+}(\U,\W)}{w_h}=\frac{a_h+ \lambda_h(\U)}{2a_h}, \\ 
{p_D^{h,-}}(x,\Delta t) &= \frac{D_h^{-}(\U,\W)}{w_h}=\frac{a_h- \lambda_h(\U)}{2a_h}.
\end{split}
\label{pD_syst_diag}
\ee
Therefore, as for the case of a relaxation approximation to a scalar conservation law, when it is possible to re-write the system in terms of characteristic variables, the GBMC algorithm keeps the grid-free property.
In this situation, starting with $n$ sets of samples
$(X^{h,0}_1,V^{h,0}_1), \ldots, (X^{h,0}_{N^h},V^{h,0}_{N^h})$, $h=1,\ldots,n$, where $V^{h,0}_i\in\{-a_h,a_h\}$ and each particle has fixed mass $m^h_i\in\{-m,m\}$, a new set of samples
$(X^h_1,V^h_1),\ldots,(X^h_{N_h},V^h_{N_h})$ is generated as presented in Algorithm \ref{alg.RWMCsyst}. 

\begin{algorithm}[h!]
\label{alg.RWMCsyst}
\caption{GBMC for general relaxation systems in characteristic form}
\begin{enumerate}
\item Compute
$X^h_i=X^{h,0}_i+V^{h,0}_i \Delta t,\quad h=1,\ldots,n, \quad i=1,\ldots,N^h$.
\item For each family of samples $h$ do the following: 
\begin{enumerate}
\item For each sample $X^h_i$, evaluate $u_i^h(t)$ using \eqref{eqch3:cs3} or as in \eqref{eq:reco}.
\item For each sample pair $(X^h_i,V^{h,0}_i)$
\begin{enumerate}
\item with probability $1-e^{-\Delta t/\e}$ do the following:
\begin{itemize}
\item[-] with probability 
${p_D^{h,+}}(X_i,t)$ defined in \eqref{pD_syst_diag} set
$V^h_i=a_h$, 
\item[-] otherwise set $V^h_i=-a_h$.
\end{itemize}
\item otherwise set $V^h_i=V_i^{h,0}$.
\end{enumerate}
\end{enumerate}
\end{enumerate}
\end{algorithm}

\begin{remark}\label{remark_GBMCmesh}
In the general case, as already mentioned, one cannot avoid the introduction of a spatial grid, since we need to reconstruct $\W$ over the grid $j$. More precisely, once this is done, for each particle $X_i$ we define the probabilities 
\be 
\begin{split}
p_D^{h,+} (X_i,\Delta t) &= \frac{|D_h^{+}(\U(X_i),\W_j)|}{|D_h^{+}(\U(X_i),\W_j)|+|D_h^{-}(\U(X_i),\W_j)|},\\
p_D^{h,-} (X_i,\Delta t) &= \frac{|D_h^{-}(\U(X_i),\W_j)|}{|D_h^{+}(\U(X_i),\W_j)|+|D_h^{-}(\U(X_i),\W_j)|},
\end{split}
\ee
and then apply the usual GBMC method, with the exception that for consistency we must introduce either particles (representing spatial derivatives) with different masses in each cell $j$ or a variable number of particles, analogous to the weighted Monte Carlo technique. 
\end{remark}

\subsection{Numerical examples for systems of conservation laws}
Here we test the two numerical methodologies, standard MC and GBMC, with two different systems of conservation laws (namely, shallow water equations and Aw-Rascle traffic model), both of which can be rewritten in terms of characteristic variables, thus allowing the gradient approach to be applied without introducing the spatial mesh. Then, we present a test case for the isentropic Euler system without the passage through the characteristic form of the model, so considering the mesh-dependent version of the GBMC method. We restrict the presentation of our results to the limiting case $\e\to 0$.

\subsubsection{Shallow water equations}
First, we consider the 1D shallow water equations with horizontal bottom topography~\cite{Toro2001}
\be
\begin{aligned}
\frac{\partial h}{\partial t} + \frac{\partial \left(h u\right)}{\partial x} &=  0,\\
\frac{\partial \left(h u\right)}{\partial t} + \frac{\partial}{\partial x}\left( \frac{gh^2}{2} + hu^2\right) &= 0,
\end{aligned}
\label{eq:SWE}
\ee
where $h$ is the water depth, $u$ is the velocity, $m=\rho u$ is the momentum, and $g$ is the gravity. 
System \eqref{eq:SWE} can be written in compact form \eqref{eq:scons}, being 
\[\U = \left(h,\, h u\right)^T, \qquad \FF(\U)=\left(h u,\, \frac{gh^2}{2} + hu^2\right)^T.\] 
Therefore, we may write the relaxation approximation \eqref{eq:reln} and then directly apply the Monte Carlo method previously discussed for the case of systems of conservation laws.

To apply the GBMC approach in mesh-less form, it is first necessary to rewrite system \eqref{eq:SWE} in terms of characteristic variables. 
Evaluating the eigenstructure of the system, whose eigenvalues result $\lambda_{1,2}=u \pm c$, where $c = \sqrt{gh}$, we can derive the following Riemann Invariants:
\be
\Gamma_{1,2} = u \pm 2c.
\label{eq.RI_SWE}
\ee
It is therefore possible to re-write the system in diagonal form, knowing that
\[\partial_t \Gamma_{1,2} + \lambda_{1,2}\,\partial_x \Gamma_{1,2} = 0,\]
so the final system reads:
\be
\begin{aligned}
\frac{\partial }{\partial t}\left( u + 2c\right)  + \left( u+c\right)\frac{\partial}{\partial x} \left(u + 2c\right) &=  0,\\
\frac{\partial  }{\partial t}\left( u - 2c\right) + \left( u-c\right)\frac{\partial}{\partial x} \left(u - 2c\right) &=  0.
\end{aligned}
\label{eq:diag_SWE}
\ee
If we define
\[
\hat\U = \begin{pmatrix} u + 2c\\ u - 2c\end{pmatrix} 
= \begin{pmatrix} \Gamma_1 \\ \Gamma_2\end{pmatrix} , \qquad
\FF'(\hat\U) = \begin{pmatrix} u+c & 0 \\ 0 & u-c \end{pmatrix}
= \begin{pmatrix} \frac{\Gamma_1 + \Gamma_2}{2} + \frac{\Gamma_1 - \Gamma_2}{4} & 0 \\ 0 & \frac{\Gamma_1 + \Gamma_2}{2} - \frac{\Gamma_1 - \Gamma_2}{4} \end{pmatrix},
\]
the system results written as
\be
\frac{\partial \hat\U}{\partial t} + \FF'(\hat\U)\frac{\partial \hat\U}{\partial x}=0.
\label{eq:syst_diagRW}
\ee
When introducing $\W=\partial \hat\U / \partial x$ and $\Z=\partial \hat\V / \partial x$, the relaxation approximation of the above system reads as \eqref{eq:reln2}, and the mesh-less GBMC algorithm can be straightforwardly applied.

\begin{figure}[!htb]
\centering
\includegraphics[width=0.42\textwidth]{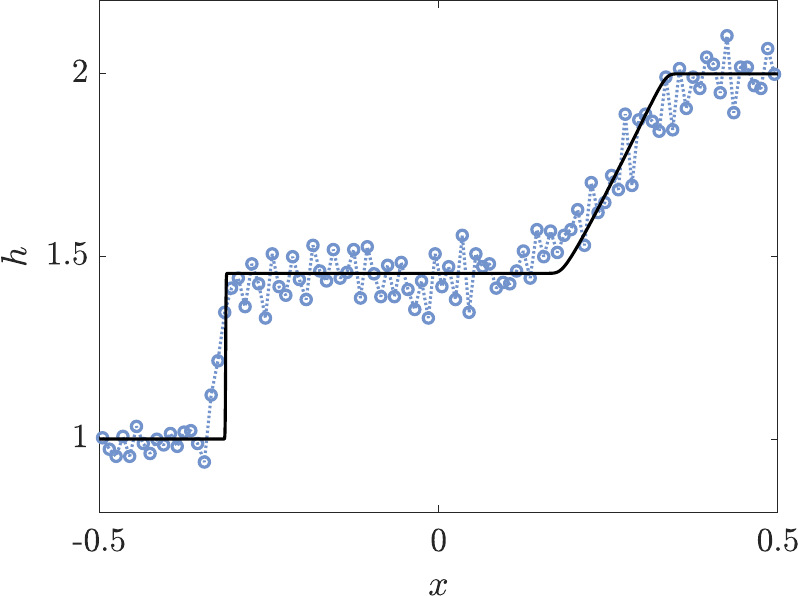}
\includegraphics[width=0.42\textwidth]{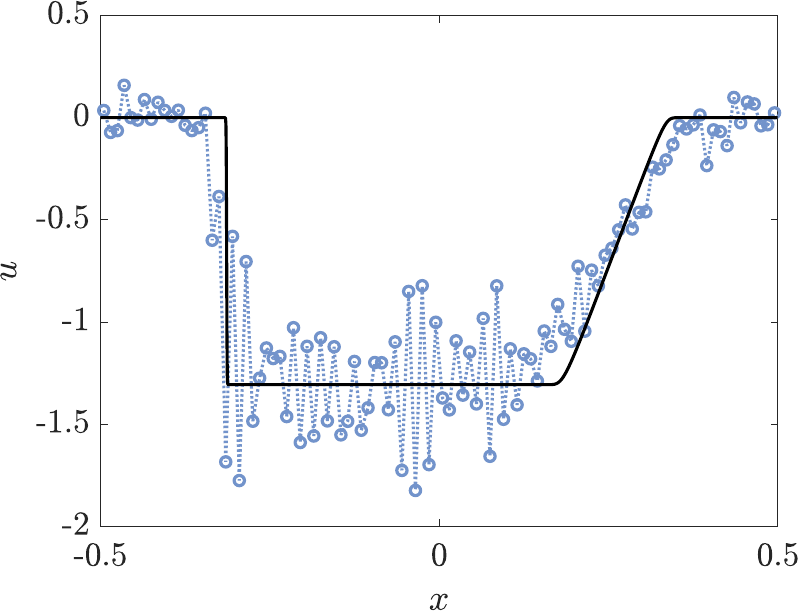}
\includegraphics[width=0.42\textwidth]{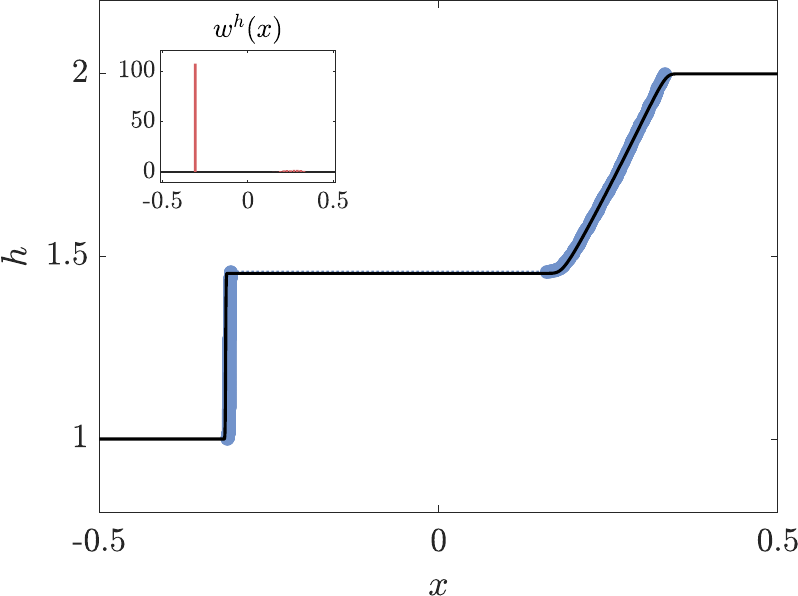}
\includegraphics[width=0.42\textwidth]{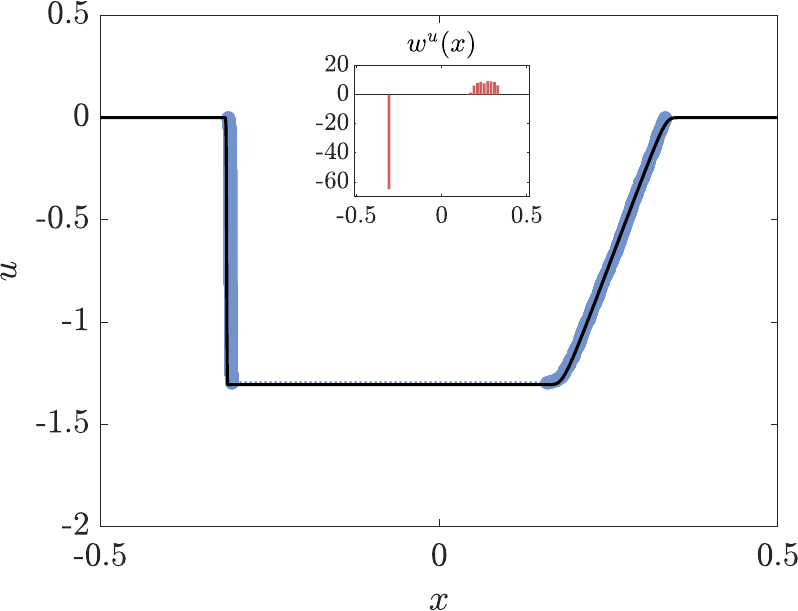}
\caption{Test 3(a), Shallow Water Equations. Solution in terms of water depth $h$ (left) and velocity $u$ (right) at $t=0.075$. 
Top: Monte Carlo with $N=100000$ particles and $M=100$ grid points, $\Delta t = 0.001$,  applying the low variance technique. 
Bottom: GBMC method with $N=2000$ particles, $\Delta t = 0.0001$. In the box, the histogram of the space derivative of the state variable is shown. 
The limit case is $\e=10^{-8}$ with $a_1=4.45$ and $a_2=5.10$. The reference solution is reported in black solid line.}
\label{SWE_RP1}
\end{figure}
\begin{figure}[!htb]
\centering
\includegraphics[width=0.42\textwidth]{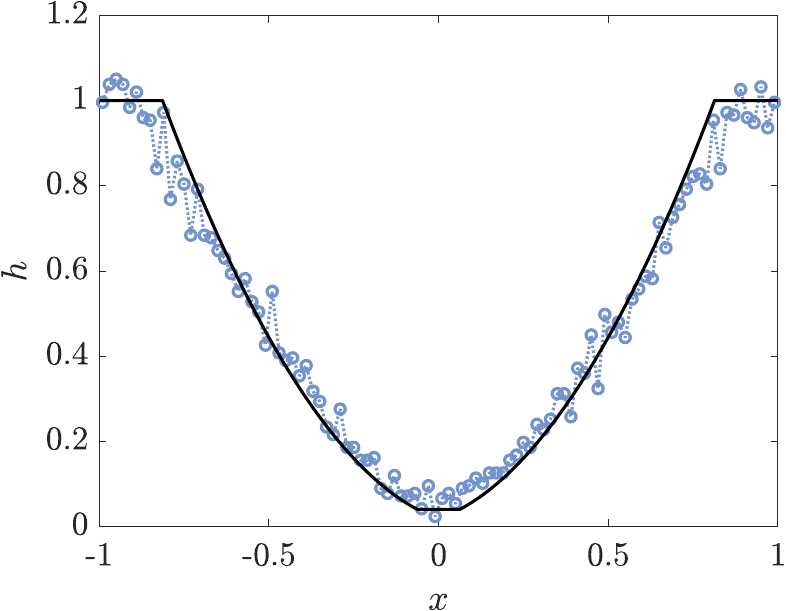}
\includegraphics[width=0.42\textwidth]{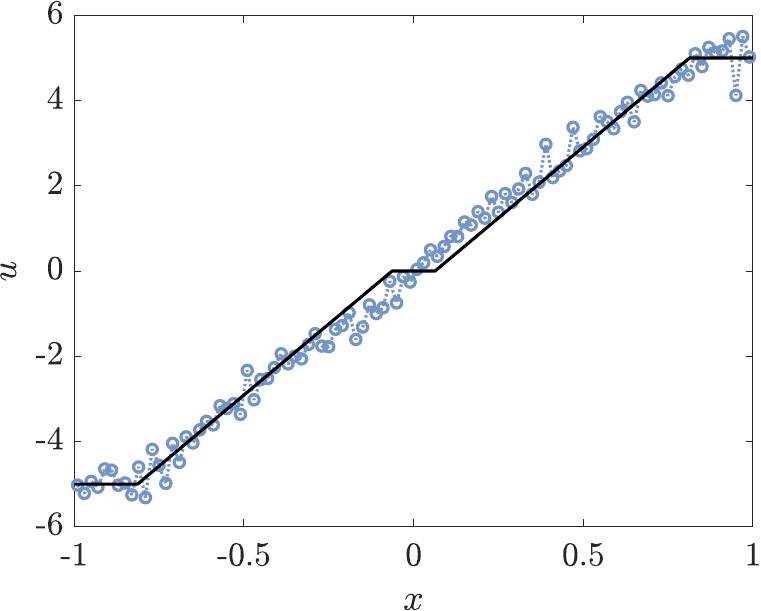}
\includegraphics[width=0.42\textwidth]{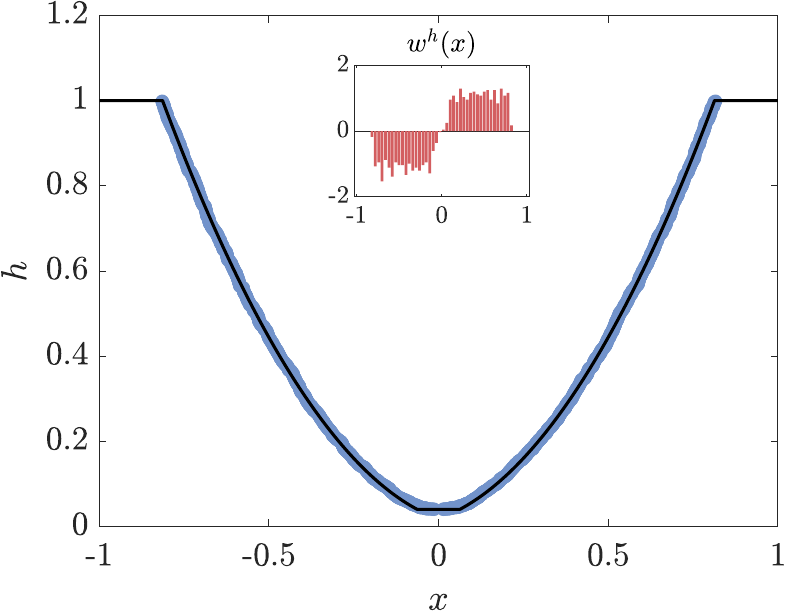}
\includegraphics[width=0.42\textwidth]{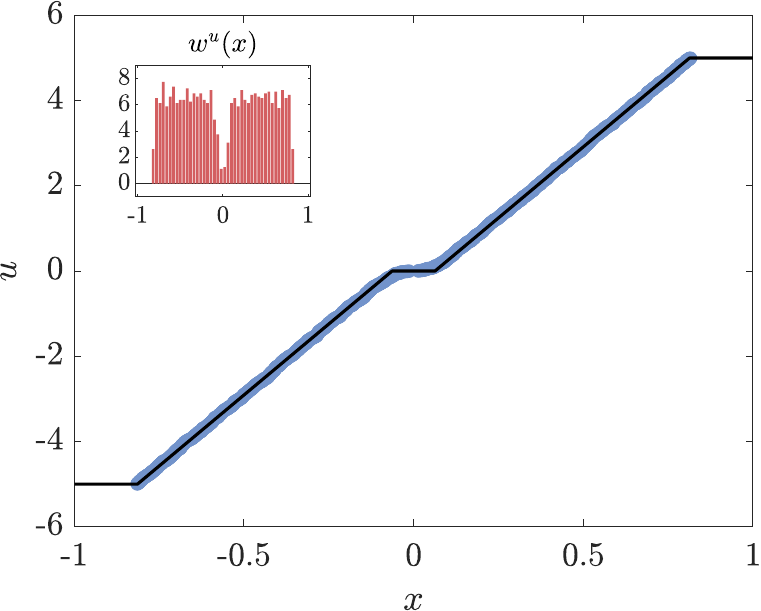}
\caption{Test 3(b), Shallow Water Equations. Solution in terms of water depth $h$ (left) and velocity $u$ (right) at $t=0.1$.
Top: Monte Carlo with $N=100000$ particles and $M=100$ grid points, $\Delta t = 0.001$, with low variance technique. 
Bottom: GBMC method with $N=2000$ particles, $\Delta t = 5\cdot 10^{-5}$. In the box, the histogram of the space derivative of the state variable is shown. 
The limit case is $\e=10^{-8}$ with $a_{1,2}=8.20$. The exact solution is reported in black solid line.}
\label{SWE_RP2}
\end{figure}

We test the Monte Carlo and GBMC methods here proposed with two Riemann Problems designed referring to~\cite{Toro2001}.
\begin{itemize}
\item\textbf{Test 3(a):} In the first case, in the domain $\Omega=[-0.5,0.5]$, we set 
\begin{equation*}
\begin{aligned}
&h(x,0) = 1, &\quad &u(x,0)=0, \qquad &\mathrm{for} \quad x<0\,,\\
&h(x,0) = 2, &\quad &u(x,0)=0, \qquad &\mathrm{for} \quad x\ge0\,.
\end{aligned}
\end{equation*}

\item\textbf{Test 3(b):} In the second case, we consider an almost dry bed solution in the domain $\Omega=[-1,1]$, imposing as initial conditions
\begin{equation*}
\begin{aligned}
&h(x,0) = 1, &\quad &u(x,0)=-5, \qquad &\mathrm{for} \quad x<0\,,\\
&h(x,0) = 1, &\quad &u(x,0)=5, \qquad &\mathrm{for} \quad x\ge0\,.
\end{aligned}
\end{equation*}
\end{itemize}

We compare the solutions obtained in terms of the variables $h$ and $u$ by both methods in Figures \ref{SWE_RP1} and \ref{SWE_RP2}. We remark here that in the GBMC approach particles are evolving in space and time following the characteristic variables $\Gamma_{1,2}$ defined in \eqref{eq.RI_SWE}. Therefore, the plots here shown are obtained considering that
\[h = \frac{c^2}{g} = \frac{(\Gamma_1 - \Gamma_2)^2}{16g} \quad \mathrm{and} \quad u = \frac{\Gamma_1 + \Gamma_2}{2}\,.\]
The augmented accuracy and the highly reduced variance of the solutions produced using the GBMC appear here even more evident when compared to results obtained with the direct Monte Carlo, even considering almost 100 times smaller amount of particles. 
Moreover, the capability of the GBMC method of better capturing the position and the sharpness of shock waves is again confirmed. Let us point out once more here how this advantage is linked to the key feature of the proposed method that allows the particles to move following the gradient of the solution. 

\subsubsection{Aw-Rascle traffic model}
In this Section we test the standard and Gradient-based Monte Carlo methods considering the Aw-Rascle traffic model~\cite{Aw}:
\be
\begin{aligned}
\frac{\partial \rho}{\partial t} + \frac{\partial \left(\rho u\right)}{\partial x} &=  0,\\
\frac{\partial}{\partial t}\left( \rho u + \rho p(\rho)\right) + \frac{\partial}{\partial x}\left( \rho u^2 + \rho u p(\rho)\right) &= 0.
\end{aligned}
\label{eq:AWR}
\ee
Here $\rho$ represents the density of cars, $u$ the velocity, and $p(\rho)=$ is a given function describing the anticipation of road conditions in front of the drivers. The system can be written in the compact form \eqref{eq:scons} defining 
\[\U = \left(\rho,\, \rho u\right)^T, \qquad \FF(\U)=\left(\rho (u + p(\rho)),\, \rho u (u + p(\rho))\right)^T,\] 
so it would be possible to directly apply the Monte Carlo method. Instead, to resort to a grid-free gradient approach, it is again necessary to express the system in terms of characteristic variables. The eigenvalues of the model are $\lambda_{1} = u$, $\lambda_{2} = u - \rho p'(\rho)$, while the Riemann Invariants result 
\be
\Gamma_1 = u + p(\rho)\, \qquad \Gamma_2 = u\,.
\label{eq.RI_AWR}
\ee
Hence, we can write the system in diagonal form as
\be
\begin{aligned}
\frac{\partial}{\partial t} \left( u + p(\rho)\right)  + u\frac{\partial}{\partial x} \left(u + p(\rho)\right) &=  0,\\
\frac{\partial u }{\partial t} + \left( u - \rho p'(\rho)\right)\frac{\partial u}{\partial x}  &=  0.
\end{aligned}
\label{eq:diag_AWR}
\ee
If we define
\[
\hat\U = \begin{pmatrix} u + p(\rho) \\ u\end{pmatrix} = \begin{pmatrix} \Gamma_1 \\ \Gamma_2\end{pmatrix} , \qquad
\FF'(\hat\U) = \begin{pmatrix} u & 0 \\ 0 & u - \rho p'(\rho) \end{pmatrix}
= \begin{pmatrix} \Gamma_2 & 0 \\ 0 & \Gamma_2 - \rho p'(\rho) \end{pmatrix},
\]
the system reads as \eqref{eq:syst_diagRW} and, following the GBMC derivation, the method can be again applied without the introduction of any spacial grid. Notice that the term $\rho p'(\rho)$ can be also expressed in term of characteristic variables, but depends on the definition of the function $p(\rho)$.

\begin{figure}[!htb]
\centering
\includegraphics[width=0.42\textwidth]{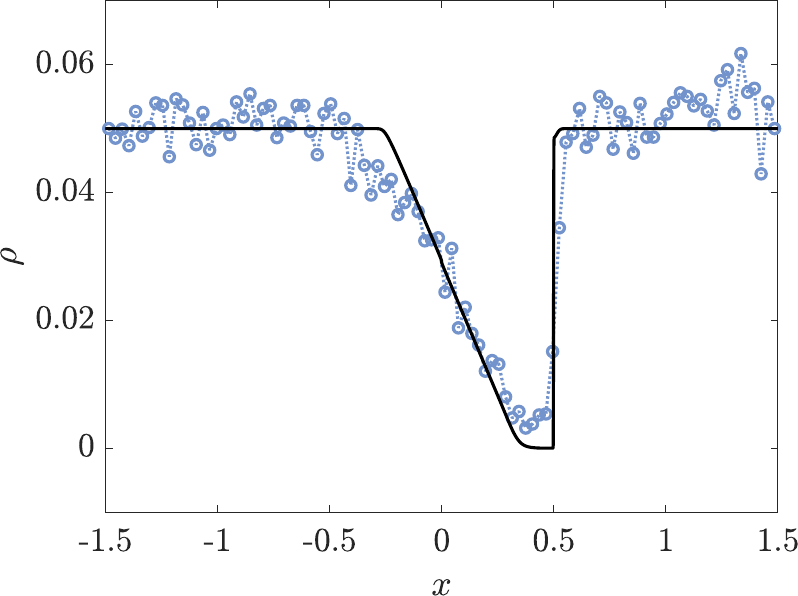}
\includegraphics[width=0.42\textwidth]{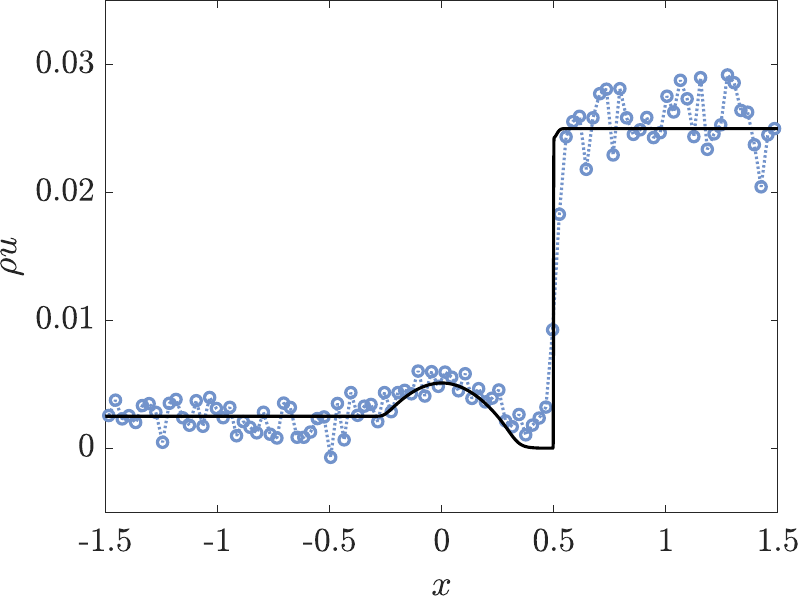}
\includegraphics[width=0.42\textwidth]{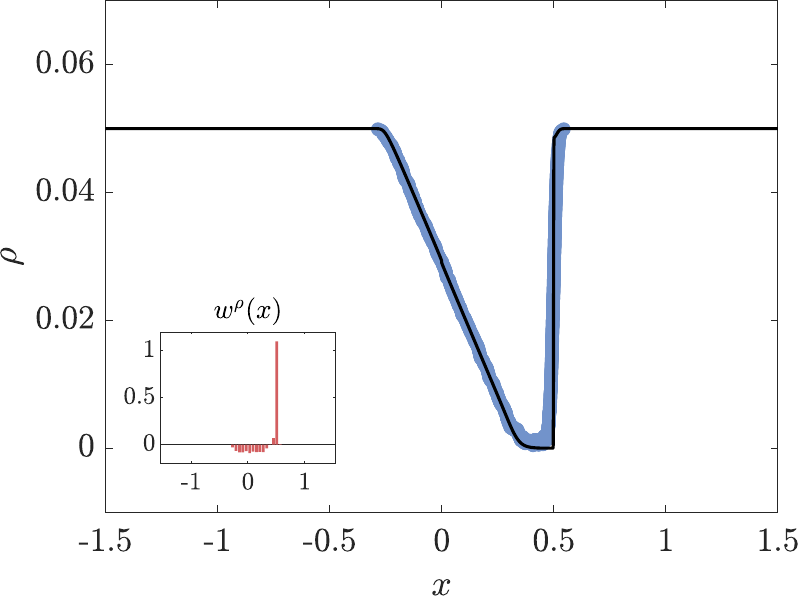}
\includegraphics[width=0.42\textwidth]{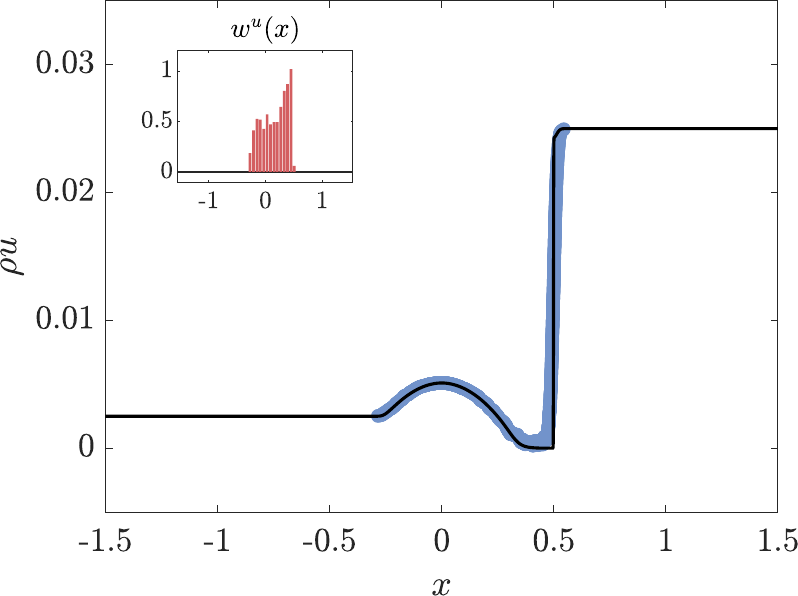}
\caption{Test 4, Aw-Rascle model. Solution in terms of density $\rho$ (left) and flux $\rho u$ (right) at $t=1.0$.
Top: Monte Carlo with $N=100000$ particles and $M=100$ grid points, $\Delta t = 0.005$, with low variance technique. 
Bottom: GBMC method with $N=2000$ particles, $\Delta t = 0.0005$. In the box, the histogram of the space derivative of the state variables $\rho$ and $u$ is shown. 
The limit case is $\e=10^{-8}$ with $a_{1,2}=0.8$. The reference solution is reported in black solid line.}
\label{AWR_RP1}
\end{figure}

We consider the Riemann Problem taken from~\cite{Aw}, which accounts for a solution with vacuum, having initial data as follows.
\begin{itemize}
\item\textbf{Test 4:} In the domain $\Omega=[-1.5,1.5]$, we have $p(\rho) = c_v \rho$, with $c_v=6$, and initial conditions
\begin{equation*}
\begin{aligned}
&\rho(x,0) = 0.05, &\quad &u(x,0)=0.05, \qquad &\mathrm{for} \quad x<0\,,\\
&\rho(x,0) = 0.05, &\quad &u(x,0)=0.5, \qquad &\mathrm{for} \quad x\ge0\,.
\end{aligned}
\end{equation*}
\end{itemize}

Notice that in this test $p'(\rho)=c_v$, hence in $\FF'(\hat\U)$ we have $F'_{22} = \Gamma_2 - \rho c_v$. 

We present the results obtained solving the problem with the standard Monte Carlo and the Gradient-based Monte Carlo in Figure \ref{AWR_RP1}, which show again an excellent behavior, even in this very challenging RP, especially for what concerns the almost absent variance in the GBMC. Notice that in the latter particles follow the dynamics of the characteristic variables $\Gamma_{1,2}$ defined in \eqref{eq.RI_AWR}. Hence, while $u =\Gamma_2$, to compute the density we need to consider that
$\rho = (\Gamma_1 - \Gamma_2)/c_v$.

\subsubsection{Isentropic Euler system}
We finally consider the following one-dimensional isentropic Euler system \cite{CJR}:
\be
\begin{aligned}
\frac{\partial \rho}{\partial t} + \frac{\partial \left(\rho u\right)}{\partial x} &=  0,\\
\frac{\partial \left(\rho u\right)}{\partial t} + \frac{\partial}{\partial x}\left( \frac1{2}\left(\rho + \rho u^2\right)\right) &= 0,
\end{aligned}
\label{eq:simple_Euler}
\ee
where $\rho$ is the gas density, $u$ is the velocity, and $m=\rho u$ is the momentum. 
System \eqref{eq:simple_Euler} can be written in the compact form \eqref{eq:scons}, being 
\[\U = \left(\rho,\, \rho u\right)^T, \qquad \FF(\U)=\left(\rho u,\, \frac1{2}\left( \rho + \rho u^2\right)\right)^T.\] 
Therefore, we may write the relaxation approximation \eqref{eq:reln} and then apply either the standard Monte Carlo approach or the Gradient method, the latter in the mesh-dependent version discussed in Remark \ref{remark_GBMCmesh}, since the system is not diagonal.

We test the resolution of the isentropic Euler equations with the two numerical methods considering two Riemann Problems having the following initial conditions, taken from \cite{CJR}.
\begin{itemize}
\item\textbf{Test 5(a):} In the first case, in the domain $\Omega=[-1,1]$, we set 
\begin{equation*}
\begin{aligned}
&\rho(x,0) = 2, &\quad &m(x,0)=1, \qquad &\mathrm{for} \quad x<0.2\,,\\
&\rho(x,0) = 1, &\quad &m(x,0)=0.13962, \qquad &\mathrm{for} \quad x\ge0.2\,.
\end{aligned}
\end{equation*}

\item\textbf{Test 5(b):} In the second case, again considering the domain $\Omega=[-1,1]$, we impose 
\begin{equation*}
\begin{aligned}
&\rho(x,0) = 1, &\quad &m(x,0)=0, \qquad &\mathrm{for} \quad x<0\,,\\
&\rho(x,0) = 0.2, &\quad &m(x,0)=0, \qquad &\mathrm{for} \quad x\ge0\,.
\end{aligned}
\end{equation*}
\end{itemize}

For both tests, we also show solutions obtained by implementing the variance reduction technique presented in Algorithm \ref{alg.VarReduction} for the standard Monte Carlo. 

\begin{figure}[!htbp]
\centering
\includegraphics[width=0.42\textwidth]{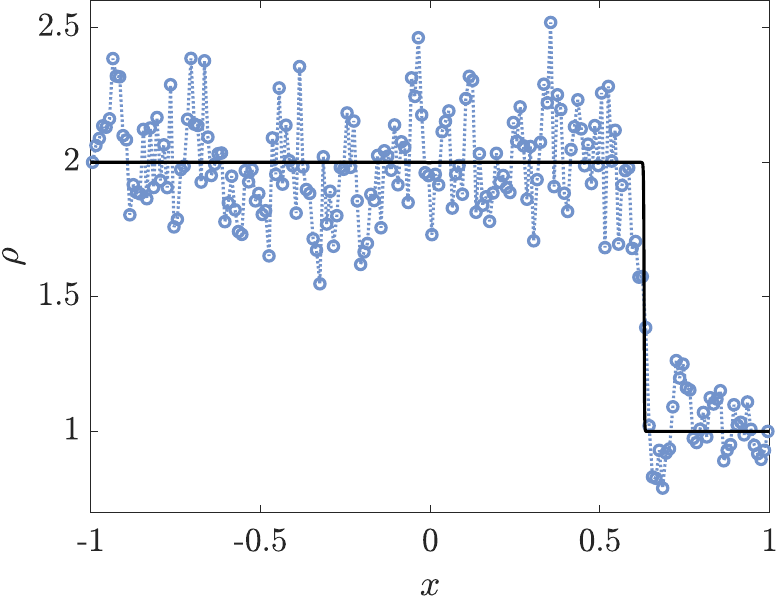}
\includegraphics[width=0.42\textwidth]{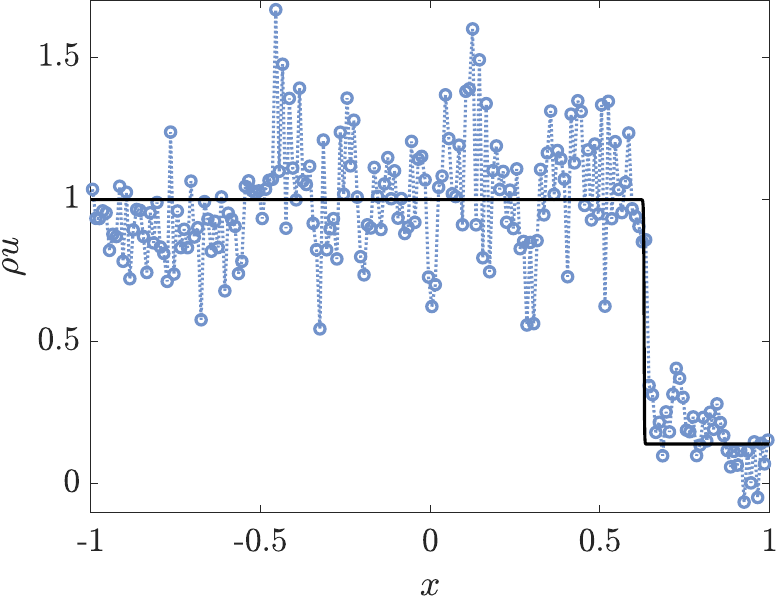}
\includegraphics[width=0.42\textwidth]{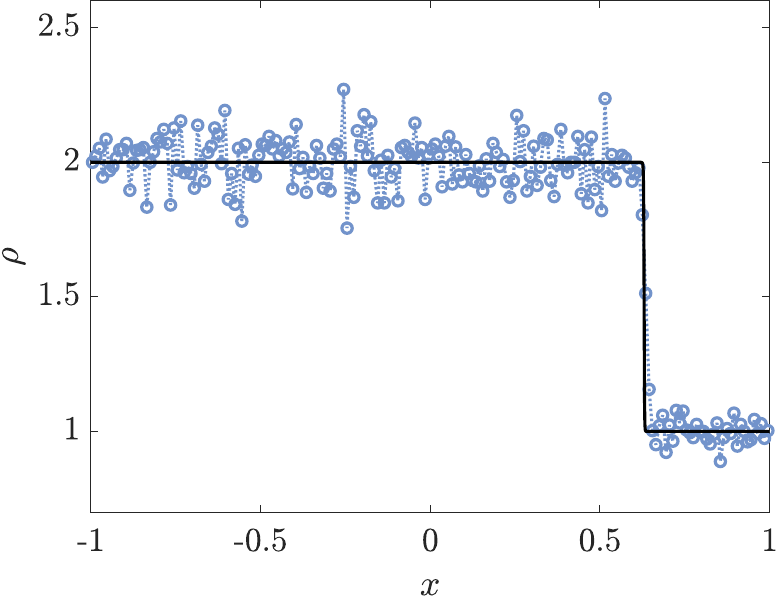}
\includegraphics[width=0.42\textwidth]{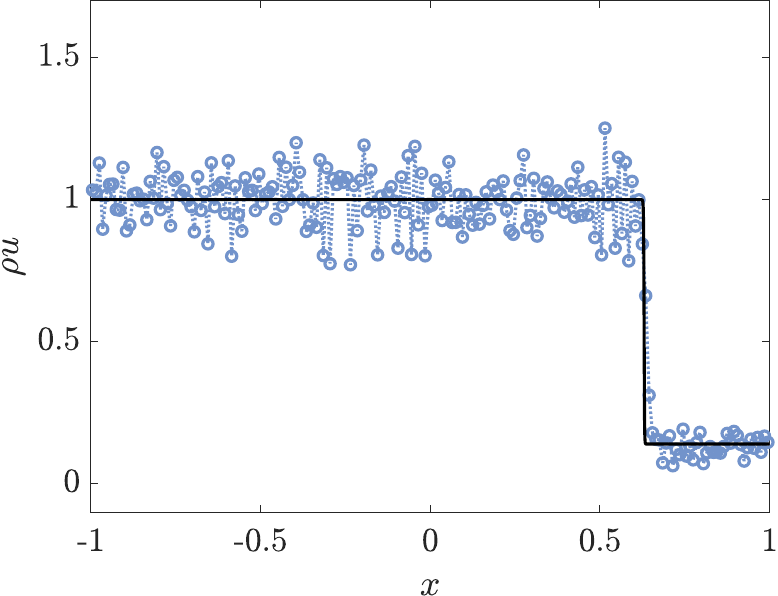}
\includegraphics[width=0.42\textwidth]{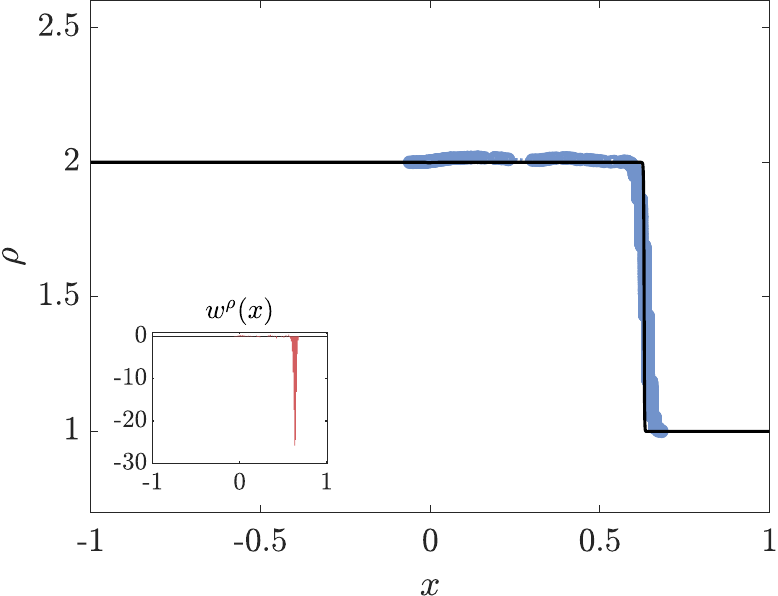}
\includegraphics[width=0.42\textwidth]{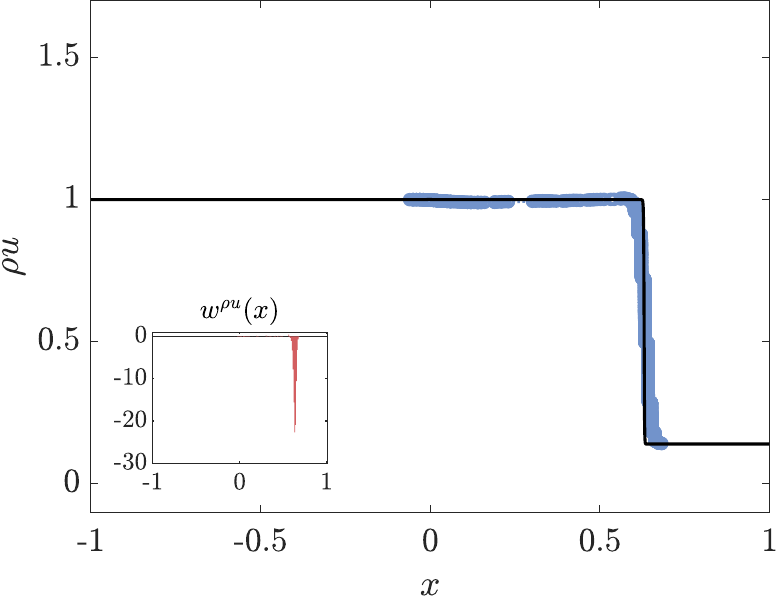}
\caption{Test 5(a), isentropic Euler system. Solution in terms of density $\rho$ (left) and flux $\rho u$ (right) at $t=0.5$.
Top: Monte Carlo with $N=200000$ particles and $M=200$ grid points, $\Delta t = 0.005$, without low variance technique. Middle: same as previous row, but with low variance technique.
Bottom: GBMC method in mesh-dependent version with $N=2000$ particles and $M=200$ grid points, $\Delta t = 0.005$. In the box, the histogram of the space derivative of the state variables  is shown. 
The limit case is $\e=10^{-8}$ with $a_{1,2}=1.0$. 
The reference solution is reported in black solid line.}
\label{Euler_RP1}
\end{figure}
\begin{figure}[!htbp]
\centering
\includegraphics[width=0.42\textwidth]{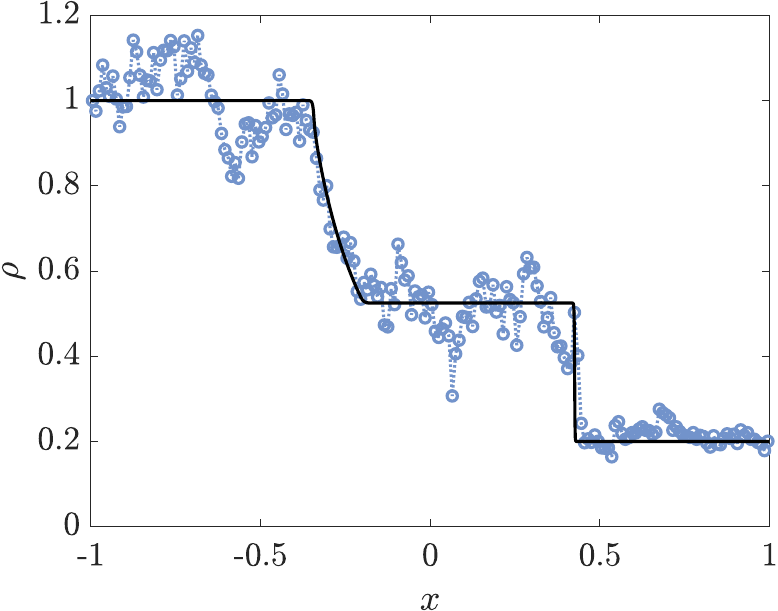}
\includegraphics[width=0.42\textwidth]{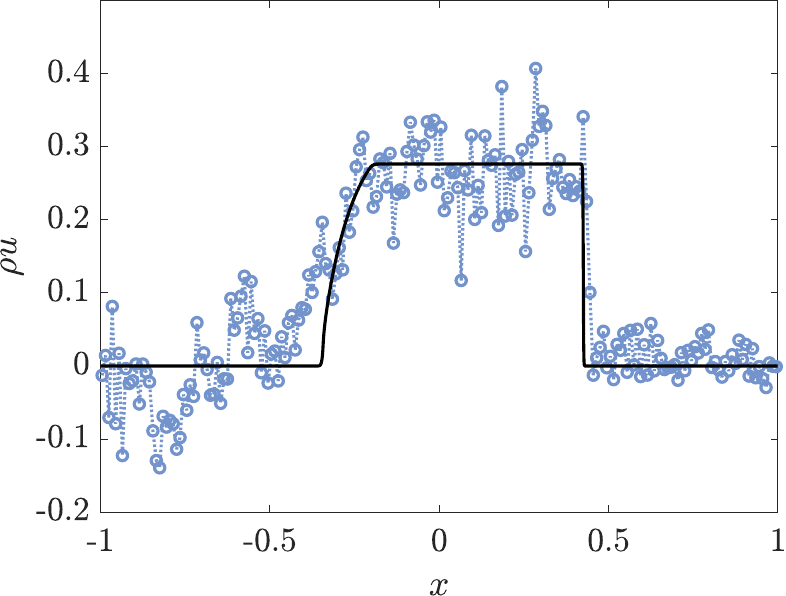}
\includegraphics[width=0.42\textwidth]{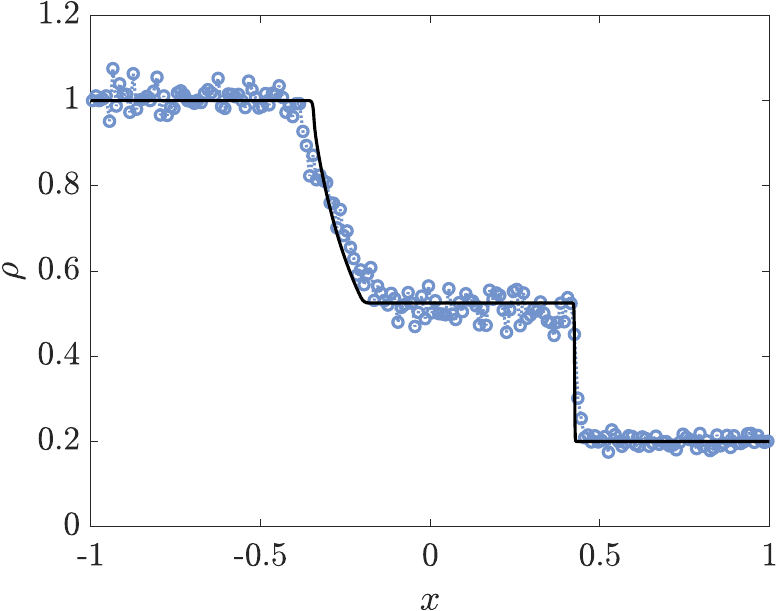}
\includegraphics[width=0.42\textwidth]{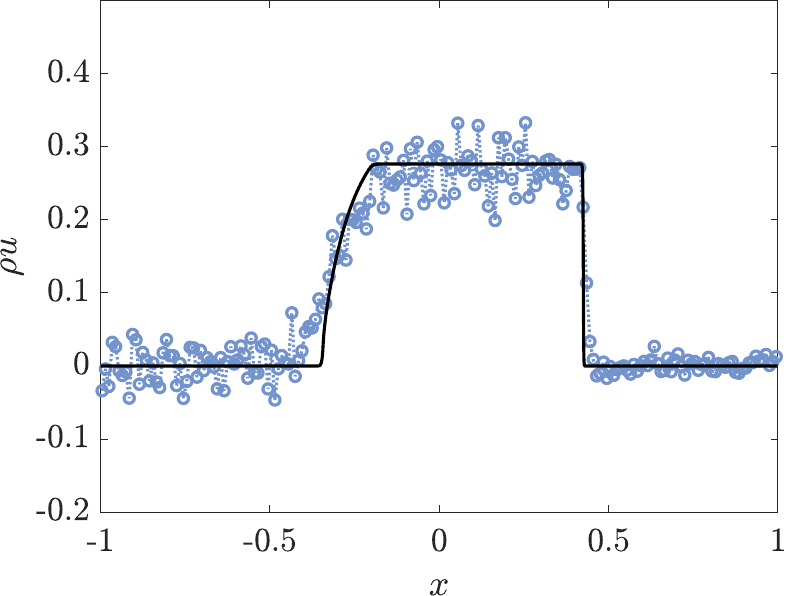}
\includegraphics[width=0.42\textwidth]{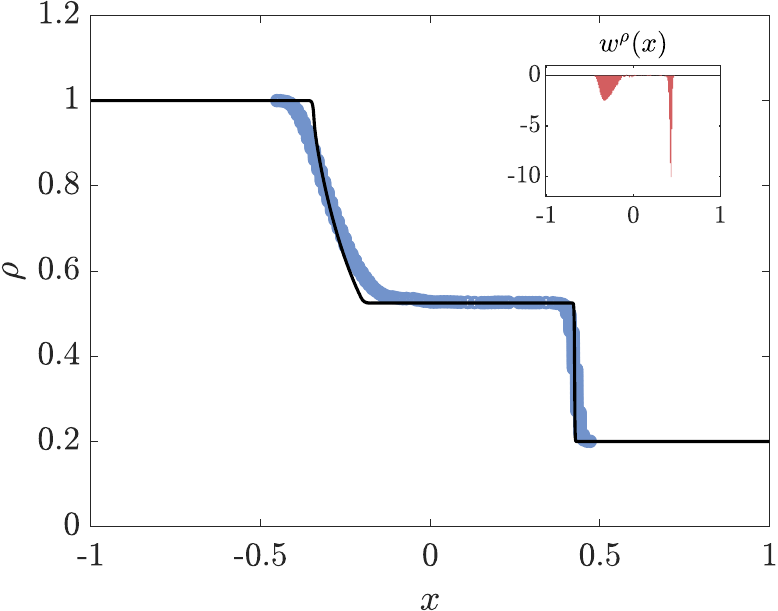}
\includegraphics[width=0.42\textwidth]{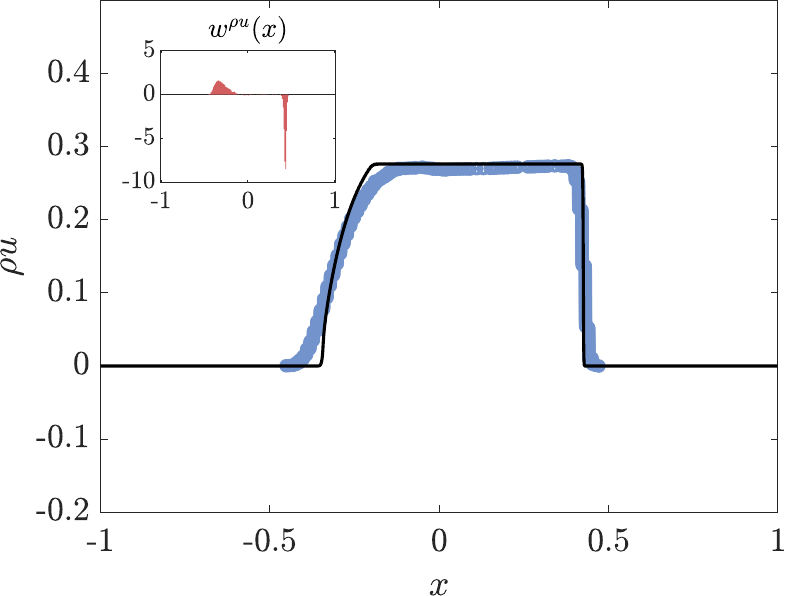}
\caption{Test 5(b), isentropic Euler system. Solution in terms of density $\rho$ (left) and flux $\rho u$ (right) at $t=0.5$.
Top: Monte Carlo with $N=200000$ particles and $M=200$ grid points, $\Delta t = 0.005$, without low variance technique. Middle: same as previous row, but with low variance technique.
Bottom: GBMC method in mesh-dependent version with $N=5000$ particles and $M=200$ grid points, $\Delta t = 0.005$. In the box, the histogram of the space derivative of the state variables  is shown. 
The limit case is $\e=10^{-8}$ with $a_{1,2}=1.0$. 
The reference solution is reported in black solid line.}
\label{Euler_RP2}
\end{figure}

Once more, looking at Figures \ref{Euler_RP1} and \ref{Euler_RP2}, the augmented accuracy and the highly reduced variance of the solutions produced with the Gradient approach appear evident when compared to results obtained with the standard Monte Carlo method, even considering a 10 times smaller amount of particles and the same spatial grid. In these figures, it is also possible to appreciate the beneficial effects of the variance reduction technique applied to the standard Monte Carlo. 
Moreover, the capability of the GBMC of better capturing the position and the sharpness of shock waves is again confirmed. Let us point out once more here that this advantage is linked to the key feature of the proposed method, which allows particles to move following the gradient of the solution.

\section{Conclusions}
Monte Carlo methods have become increasingly important in scientific computing due to their ability to handle complex systems and quantify uncertainty. Despite this, their systematic use in solving partial differential equations is still limited compared to deterministic approaches, which in many cases provide greater flexibility and accuracy. In this paper, we attempt to take a step forward in the design of Monte Carlo methods for PDEs by analyzing their systematic use for relaxation approximations to systems of hyperbolic conservation laws \cite{JinXin}. On the one hand, we extend Monte Carlo techniques of direct simulation inspired by kinetic theory to systems of hyperbolic conservation laws. On the other hand, we consider a different approach based on the use of the spatial derivative of the solution, which was developed earlier for reaction-diffusion equations \cite{SP86}, and refer to as Gradient-based Monte Carlo (GBMC). The latter method has shown great potential due to its ability to concentrate particles where the solution has large derivatives and has a grid-free structure. In the presented test cases, the GBMC method proved to be significantly more accurate than the Monte Carlo method.

Lest us notice that, by combining the techniques here presented to previous results for reaction-diffusion problems it is possible to deal with general systems of PDEs involving convection-diffusion-reaction terms. In this paper, we have limited our analysis to the one-dimensional case. In the future, our aim is to extend the GBMC method to the multidimensional case through a component-wise approach and to more general systems. Another interesting direction is to extend the GBMC approach to kinetic equations, including the challenging case of the Boltzmann equation. We leave further investigations on these topics to future research.
\section*{Acknowledgments}
This work has been written within the activities of GNCS group of INdAM (Italian National Institute of High Mathematics), whose support is acknowledged. It has also been partially supported by ICSC -- Centro Nazionale di Ricerca in High Performance Computing, Big Data and Quantum Computing, funded by European Union -- NextGenerationEU and by MIUR-PRIN Project 2017, No. 2017KKJP4X ``Innovative numerical methods for evolutionary partial differential equations and applications". G.B. would like to thank the Courant Institute, New York University, for the kind hospitality during her visiting period, and the support of the University of Ferrara with the mobility fund under ``Bando Giovani Ricercatori 2019''.


\section*{Appendix}
\appendix

\section{Reconstruction error estimates}
\label{appendix}
In this section we will detail how the numerical solution is recovered starting from the particles and illustrate how this will affect the reconstruction error both for MC and GBMC. The control of the distance between some reconstruction of the empirical measure and its true distribution is, of course, a long standing problem, central both in probability, statistics and computer science. While many distances can be used to consider the problem, here we will focus on the $L^p$ distance given its central use in error estimation for conservation laws. However, other distances, such as the Wasserstein distance, are quite natural for quantifying particle approximations of PDEs^^>\cite{Fournier15}. In the sequel, to simplify notations, we consider only one family of particles, so ignore their direction of motion, and restrict to nonnegative probability densities $u(x,t) \geq 0$.

\subsection{Monte Carlo} 
\label{appendix_MC}
For the standard Monte Carlo method, starting from $N$ samples $X_1$, \ldots, $X_N$ at time $t$ i.i.d. as $u(x,t)$ we compute the empirical density function
\begin{equation}
u_N(x,t)=\frac1{N}\sum_{k=1}^N \delta(x-X_k(t)),
\label{eq:emp}
\end{equation}
where $\delta(\cdot)$ is the Dirac delta function.

The above function needs to be regularized for numerical purposes. To this aim, let us introduce the mesh width $\Delta x>0$ and denote by $S_{\Delta x}\ge 0$ a smoothing function such that 
\[
\int_{\Omega} S_{\Delta x}(x)dx = 1.
\]
Then, we consider the approximation of the empirical density \eqref{eq:emp} obtained by a kernel density estimator in the form
\begin{equation}
\label{eq:f_RN}
u_{N,\Delta x}(x,t) = \dfrac{1}{N} \sum_{k = 1}^N S_{\Delta x}(x-X_k(t)). 
\end{equation}

In the simplest case, we have the rectangular kernel $S_{\Delta x}(x)=\chi(|x|\leq \Delta x/2)/\Delta x$, where $\chi(\cdot)$ is the indicator function, corresponding to the histogram reconstruction that has been used on a fixed grid of size $\Delta x$ in the numerical results. We refer to^^>\cite{Tsybaokv} for examples of other kernels.

%

By considering only the error due to the reconstruction  \eqref{eq:f_RN} on the computational domain $\Omega\subseteq \mathbb R$, this can be estimated from triangle inequality:
\be
\begin{split}
\left\| u(\cdot,t)-u_{N,\Delta x}(\cdot,t)\right\|_{L^p(\Omega,L^2(\Omega))} & \leq 
\| u(\cdot,t)-u_{\Delta x}(\cdot,t)\|_{L^p(\Omega)}\\
&+\| u_{\Delta x}(\cdot,t)-u_{N,\Delta x}(\cdot,t)\|_{L^p(\Omega,L^2(\Omega))},
\end{split}
\label{est1}
\ee
where
\[
u_{\Delta x}(x,t)=\int_{\Omega} S_{\Delta x}(x-y)u(y,t)\,dy
\]
and we defined
\begin{equation*}
\|g\|_{L^p(\Omega,L^2(\Omega))}=\|\mathbb {E}\left[g^2\right]^{1/2}\|_{L^p(\Omega)},
\end{equation*}
with $\mathbb{E}[\cdot]$ denoting the expectation with respect to the random variables $x_1,\ldots,x_N$ i.i.d. as $u(x,t)$ and $\|\cdot\|_{L^p(\Omega)}$ is the $L^p$ norm in $\Omega$.

Now, for the second term on the r.h.s. of \eqref{est1} we have the following lemma.
\begin{lemma} The root mean squared error satisfies
\[
\mathbb E\left[(u_{\Delta x}(x,t)-u_{N,\Delta x}(x,t))^2\right]^{1/2} =
\frac{\sigma_S(x,t)}{N^{1/2}}
\]
where 
\begin{equation}
\sigma^2_S(t,x)=\int_{\Omega}\left(S_{\Delta x}(x-y)-u_{\Delta x}(x,t)\right)^2 u(y,t)\,dy.
\label{eq:sigmas}
\end{equation}
\label{lem1}
\end{lemma} 
The proof follows by classical arguments on the convergence of the root mean squared error\cite{Caflisch98,Trimborn}. 

Next, we assume that the first term satisfies
\[
\| u(\cdot,t)-u_{\Delta x}(\cdot,t)\|_{L^p(\Omega)} \leq C_q (\Delta x)^q,
\]
according to the order of accuracy $q\geq 1$ used in the reconstruction. For example, in the case of  histogram reconstruction through a rectangular kernel we have $q=2$. In fact, from the midpoint rule we can compute
\be
\left|u(x,t)-u_{\Delta x}(x,t)\right| = \left|u(x,t)-\frac1{\Delta x}\int_{x-\Delta x/2}^{x+\Delta x/2} u(y,t)\,dy\right| 
\leq \frac1{24}\left|\frac{\partial^2 u}{\partial x^2}(\xi,t)\right|\Delta x^2,
\label{histe}
\ee
where $\xi$ depends on $x$. 

Thus, we have the following result.
\begin{theorem} For a sufficiently smooth function $u(x,t)$ the error introduced by the reconstruction function \eqref{eq:f_RN} satisfies
\begin{equation*}
\left\| u(\cdot,t)-u_{N,\Delta x}(\cdot,t)\right\|_{L^p(\Omega,L^2(\Omega))} \leq \frac{\|\sigma_S\|_{L^p(\Omega)}}{N^{1/2}} + C_q (\Delta x)^q,
\end{equation*}
where $C_{q}$ depends on the $q$ derivative of $u(x,t)$ and the domain $\Omega$, and $\sigma_{S}^2$ is defined in \eqref{eq:sigmas}.
\label{th:1}
\end{theorem}
However, more generally, the error is affected by the numerical solution of the PDE, which in our case is first order in time, due to the time splitting algorithm, and first order in space, due to the piecewise constant reconstruction used in the relaxation phase. Roughly speaking, we can assume that these errors, as a result of a CFL-type condition, are of the same order and to have a local error estimate of the type
\be
\left|u(x,t)-\tilde u_{\Delta x}(x,t)\right| \leq M_1(x)\Delta x,
\label{histe3}
\ee
where $\tilde u_{\Delta x}(x,t)=\mathbb{E}[\tilde u_{N,\Delta x}(x,t)]$, with $\tilde u_{N,\Delta x}(x,t)$ the reconstructed MC solution of the PDE. In \eqref{histe3} the quantity $M_1(x)$ depends on the PDE under consideration and the first order derivative of the solution. Note that, since the first order spatial error in solving the PDE dominates the reconstruction error in \eqref{histe}, we can ignore in \eqref{histe3} the second order contributions.

Now, we can estimate the MC solution of the PDE from 
\be
\begin{split}
\left\| u(\cdot,t)-\tilde u_{N,\Delta x}(\cdot,t)\right\|_{L^p(\Omega,L^2(\Omega))} & \leq 
\| u(\cdot,t)-\tilde u_{\Delta x}(\cdot,t)\|_{L^p(\Omega)}\\
&+\| \tilde u_{\Delta x}(\cdot,t)-\tilde u_{N,\Delta x}(\cdot,t)\|_{L^p(\Omega,L^2(\Omega))},
\end{split}
\label{est2}
\ee
where the second term can be estimated as in Lemma \ref{lem1}:
\[
\mathbb E\left[(\tilde u_{\Delta x}(x,t)-\tilde u_{N,\Delta x}(x,t))^2\right]^{1/2} =
\frac{\tilde\sigma(x,t)}{N^{1/2}}
\]
with 
\begin{equation}
\begin{split}
\tilde\sigma^2(t,x)&=\mathbb{E}[(\tilde u_{N,\Delta x}(x,t))^2]-\mathbb{E}[\tilde u_{N,\Delta x}(x,t)]^2\\
&\leq \frac{\tilde u_{\Delta x}(x,t)}{\Delta x}-(\tilde u_{\Delta x}(x,t))^2\leq \frac{\tilde u_{\Delta x}(x,t) }{\Delta x}\leq \frac{u(x,t)}{\Delta x}+M_1(x),
\end{split}
\label{eq:sigmas2}
\end{equation}
and the last inequality follows from assumption \eqref{histe3}.

By collecting the above results we can state the following theorem.
\begin{theorem} Let us denote by $\tilde u_{N,\Delta x}(\cdot,t)$ the reconstructed MC solution such that it satisfies \eqref{histe3}, then  
\be
\left\| u(\cdot,t)-\tilde u_{N,\Delta x}(\cdot,t)\right\|_{L^p(\Omega,L^2(\Omega))} \leq \frac{\|u\|^{1/2}_{L^{p/2}(\Omega)}}{(\Delta x N)^{1/2}} + \frac{C_2}{N^{1/2}}+C_1 \Delta x
\label{eq:emc}
\ee
where the constants $C_1$, $C_2$ depend on the first order derivative of $u(x,t)$ and the domain $\Omega$.
\end{theorem}
It is thus clear that we must take $\Delta x$ small to make the mesh bias small, but also large enough to ensure that the variance is also small. If we want to minimize the error with respect to $\Delta x$ we should
take
\be
\Delta x = \left(\frac{\|u\|_{L^{p/2}(\Omega)}}{4C_1^2 N}\right)^{1/3}.
\label{eq:dxopt}
\ee
This will ensure an optimal error decay $O(N^{-1/3})$.  Therefore, the optimal mesh will scale as $\Delta x \approx N^{-1/3}$ for a convergence rate of $O(N^{-1/3})$, where the precise value of $C_1$ depends on the particular PDE under consideration and may be difficult to estimate in practice.    
%

\subsection{Gradient-based Monte Carlo} 
Assume, again for simplicity of treatment, that $u(x,t)$ is nondecreasing monotone so that $w(x,t)=\frac{\partial u(x,t)}{\partial x} \geq 0$. Therefore, without loss of generalization, we can consider $w(x,t)$ as a probability density and $u(x,t)$ its cumulative distribution function. 

Given now $N$ samples $X_1$, \ldots, $X_N$ at time $t$ i.i.d. as $w(x,t)$, we compute the empirical cumulative distribution function as
\begin{equation}
u_N(x,t)=\frac1{N}\sum_{k=1}^N H(x-X_k(t)),
\label{eq:emp2}
\end{equation}
where $H(\cdot)$ is the Heaviside step function.

Of course, \eqref{eq:emp2} differs from \eqref{eq:emp} although we used the same notations for the two reconstructions. Note that, in contrast to \eqref{eq:emp}, the empirical cumulative distribution does not need any further approximation or the introduction of a mesh but can be used directly in the form \eqref{eq:emp2}.  

We recall some classical results for \eqref{eq:emp2} (see^^>\cite{vv98} for more details).
\begin{lemma}
For any given $x\in\Omega$ we have that $Nu_N(x,t)$ has a binomial distribution with parameters $N$ and success probability $u(x,t)$. Therefore
\be
\mathbb{E}[u_N(x,t)] = u(x,t),\qquad {\rm \mathbb{V}ar}[u_N(x,t)]=\frac1{N}u_N(x,t)(1-u_N(x,t)),
\ee
where ${\rm \mathbb{V}ar}[\cdot]$ denotes the variance with respect to the random variables $x_1,\ldots,x_N$ i.i.d. as $w(x, t)$, and $u_N(x,t)$ converges to $u(x,t)$ almost surely.
\end{lemma}
In fact, $u_N(x,t)$ is the sum of $N$ independent Bernoulli random variables, therefore $N u_N(x,t)$ is a binomial random variable. Thus, the mean and variance characterization in the proposition follows easily. The latter statement is a consequence of Hoeffding’s inequality, which implies that for any $\e>0$ in probability we have
\be
{\mathcal P}\left(|u_N(x,t)-u(x,t)|\geq \e \right) \leq 2e^{-2N\e^2}.
\ee 
The above pointwise convergence can be made uniform by the fundamental Glivenko-Cantelli lemma.
\begin{lemma}[Glivenko-Cantelli] The empirical distribution $u_N(x,t)$ converges uniformly to $u(x,t)$, namely as $N\to\infty$ we have
\be
{\rm sup}_{x\in\Omega}|u_N(x,t)-u(x,t)|\overset{a.s.}{\to} 0
\ee
where the superscript a.s. denotes convergence almost surely.
\end{lemma}
A similar estimate as Hoeffding’s inequality holds true also in this latter case.

Now, let's consider the problem of estimating the convergence rate in $L^p$ spaces. Since we have
\[
u(x,t)=\int_{\Omega} H(x-y)w(y)\,dy,\qquad u_N(x,t)=\int_{\Omega} H(x-y)w_N(y)\,dy,
\]
the numerical error of the empirical cumulative distribution \eqref{eq:emp2}, can now be estimated similarly to Lemma \ref{lem1} from
\be
\begin{split}
\mathbb{E}[(u(x,t)-u_{N}(x,t))^2]^{1/2} = \frac{\sigma_H(x,t)}{N^{1/2}}
\end{split}
\label{est3}
\ee
where
\begin{equation}
\sigma^2_H(x,t)=\int_{\Omega}\left(H(x-y)-u(x,t)\right)^2 w(y,t)\,dy.
\label{eq:sigmah}
\end{equation}
Thus, we have the following theorem.
\begin{theorem} The empirical cumulative distribution \eqref{eq:emp2} satisfies 
\be
\left\| u(\cdot,t)-u_{N}(\cdot,t)\right\|_{L^p(\Omega,L^2(\Omega))} \leq \frac{\|u\|^{1/2}_{L^{p/2}(\Omega)}}{N^{1/2}}.
\label{eq:egbmc}
\ee
\end{theorem}
The proof follows immediately by observing that
\[
\sigma^2_H(x,t) = u(x,t)(1-u(x,t)) \leq u(x,t).
\]
If we now consider the error introduced by the GBMC solution of the PDE, in the case when this is affected only by the time error of the splitting since no mesh in space is needed, we can assume  
\be
\left\|u(\cdot,t)-\tilde u(\cdot,t)\right\|_{L^p(\Omega)} \leq \tilde C_1\Delta t,
\label{histe4}
\ee
where $\tilde u(\cdot,t)=\mathbb{E}[\tilde u_{N}(\cdot,t)]$, $\tilde u_{N}(\cdot,t)$ is the numerical GBMC solution of the PDE and $\tilde C_1$ depends on the first order time derivative of the solution and the domain $\Omega$.

This leads immediately to the following result.
\begin{theorem} Let us denote by $\tilde u_{N}(\cdot,t)$ the GBMC solution such that it satisfies \eqref{histe4}, then
\be
\left\| u(\cdot,t)-\tilde u_{N}(\cdot,t)\right\|_{L^p(\Omega,L^2(\Omega))} \leq \frac{\|u\|^{1/2}_{L^{p/2}(\Omega)}}{N^{1/2}}+ {\tilde C_1}{\Delta t}.
\label{eq:egbmc2}
\ee
\end{theorem}
From \eqref{eq:egbmc2} it is clear that the time step should be taken as $\Delta t \approx N^{-1/2}$ to optimize the error in the GBMC solution and achieve an optimal convergence rate of $O(N^{-1/2})$. Note that, small values of $\Delta t$ will have a moderate impact on the computational cost of GBMC since less particles interact and no extra cost of the reconstruction is needed.

By comparing \eqref{eq:egbmc2} with \eqref{eq:emc} one can see that the error for the GBMC method is always smaller than the MC error. More precisely, the two errors have the same decay for a fixed (non optimal) $\Delta x$ in the MC method, whereas for an optimal $\Delta x$ given by \eqref{eq:dxopt} the GBMC method has a faster convergence rate as $O(N^{-1/2})$ against $O(N^{-1/3})$ of the MC method.

\end{document}